\documentclass[12pt]{amsart}
\usepackage{amscd,amsmath,amssymb,amsfonts}
\usepackage[cmtip, all]{xy}
\begin{document}
\newtheorem{thm}{Theorem}[subsection]
\newtheorem{prop}[thm]{Proposition}
\newtheorem{lem}[thm]{Lemma}
\newtheorem{cor}[thm]{Corollary}
\newtheorem{conj}[thm]{Conjecture}
\theoremstyle{definition}
\newtheorem{defin}[thm]{Definition}
\theoremstyle{remark}
\newtheorem{remk}[thm]{Remark}
\newtheorem{remks}[thm]{Remarks}
\newtheorem{exm}[thm]{Example}
\newtheorem{exms}[thm]{Examples}
\newtheorem{notat}[thm]{Notation}
\numberwithin{equation}{subsection}

\newenvironment{enum}{\makeatletter
\renewcommand{\theenumi}{\roman{enumi}}
\renewcommand{\labelenumi}{(\theenumi)}
\makeatother
\begin{enumerate}}%
{\end{enumerate}}

\newenvironment{enumera}{\makeatletter
\renewcommand{\theenumi}{\alph{enumi}}
\renewcommand{\labelenumi}{(\theenumi)}
\makeatother
\begin{enumerate}}%
{\end{enumerate}}

\newenvironment{Def}{\begin{defin}}%
{\hfill$\square$\end{defin}}

\newenvironment{rem}{\begin{remk}}%
{\hfill$\square$\end{remk}}

\newenvironment{rems}{\begin{remks}}%
{\hfill$\square$\end{remks}}

\newenvironment{ex}{\begin{exm}}%
{\hfill$\square$\end{exm}}

\newenvironment{exs}{\begin{exms}}%
{\hfill$\square$\end{exms}}

\newenvironment{note}{\begin{notat}}%
{\hfill$\square$\end{notat}}

\newcommand{\thmref}{Theorem~\ref}
\newcommand{\propref}{Proposition~\ref}
\newcommand{\corref}{Corollary~\ref}
\newcommand{\defref}{Definition~\ref}
\newcommand{\lemref}{Lemma~\ref}
\newcommand{\exref}{Example~\ref}
\newcommand{\remref}{Remark~\ref}
\newcommand{\chref}{Chapter~\ref}
\newcommand{\apref}{Appendix~\ref}
\newcommand{\secref}{Section~\ref}
\newcommand{\conjref}{Conjecture~\ref}

\newcommand{\sA}{{\mathcal A}}
\newcommand{\sB}{{\mathcal B}}
\newcommand{\sC}{{\mathcal C}}
\newcommand{\sD}{{\mathcal D}}
\newcommand{\sE}{{\mathcal E}}
\newcommand{\sF}{{\mathcal F}}
\newcommand{\sG}{{\mathcal G}}
\newcommand{\sH}{{\mathcal H}}
\newcommand{\sI}{{\mathcal I}}
\newcommand{\sJ}{{\mathcal J}}
\newcommand{\sK}{{\mathcal K}}
\newcommand{\sL}{{\mathcal L}}
\newcommand{\sM}{{\mathcal M}}
\newcommand{\sN}{{\mathcal N}}
\newcommand{\sO}{{\mathcal O}}
\newcommand{\sP}{{\mathcal P}}
\newcommand{\sQ}{{\mathcal Q}}
\newcommand{\sR}{{\mathcal R}}
\newcommand{\sS}{{\mathcal S}}
\newcommand{\sT}{{\mathcal T}}
\newcommand{\sU}{{\mathcal U}}
\newcommand{\sV}{{\mathcal V}}
\newcommand{\sW}{{\mathcal W}}
\newcommand{\sX}{{\mathcal X}}
\newcommand{\sY}{{\mathcal Y}}
\newcommand{\sZ}{{\mathcal Z}}
\newcommand{\A}{{\mathbb A}}
\newcommand{\B}{{\mathbb B}}
\newcommand{\C}{{\mathbb C}}
\newcommand{\D}{{\mathbb D}}
\newcommand{\E}{{\mathbb E}}
\newcommand{\F}{{\mathbb F}}
\newcommand{\G}{{\mathbb G}}
\renewcommand{\H}{{\mathbb H}}
\newcommand{\I}{{\mathbb I}}
\newcommand{\J}{{\mathbb J}}
\newcommand{\M}{{\mathbb M}}
\newcommand{\N}{{\mathbb N}}
\renewcommand{\P}{{\mathbb P}}
\newcommand{\Q}{{\mathbb Q}}
\newcommand{\R}{{\mathbb R}}
\newcommand{\T}{{\mathbb T}}
\newcommand{\U}{{\mathbb U}}
\newcommand{\V}{{\mathbb V}}
\newcommand{\W}{{\mathbb W}}
\newcommand{\X}{{\mathbb X}}
\newcommand{\Y}{{\mathbb Y}}
\newcommand{\Z}{{\mathbb Z}}
\renewcommand{\1}{{\mathbb{S}}}

\newcommand{\fm}{{\mathfrak m}}
\newcommand{\an}{{\rm an}}
\newcommand{\alg}{{\rm alg}}
\newcommand{\cl}{{\rm cl}}
\newcommand{\Alb}{{\rm Alb}}
\newcommand{\CH}{{\rm CH}}
\newcommand{\mc}{\mathcal}
\newcommand{\mb}{\mathbb}
\newcommand{\surj}{\twoheadrightarrow}
\newcommand{\red}{{\rm red}}
\newcommand{\codim}{{\rm codim}}
\newcommand{\rank}{{\rm rank}}
\newcommand{\Pic}{{\rm Pic}}
\newcommand{\Div}{{\rm Div}}
\newcommand{\Hom}{{\rm Hom}}
\newcommand{\im}{{\rm im}}
\newcommand{\Spec}{{\rm Spec \,}}
\newcommand{\sing}{{\rm sing}}
\newcommand{\Char}{{\rm char}}
\newcommand{\Tr}{{\rm Tr}}
\newcommand{\Gal}{{\rm Gal}}
\newcommand{\Min}{{\rm Min \ }}
\newcommand{\Max}{{\rm Max \ }}
\newcommand{\supp}{{\rm supp}\,}
\newcommand{\0}{\emptyset}
\newcommand{\sHom}{{\mathcal{H}{om}}}

\newcommand{\Nm}{{\operatorname{Nm}}}
\newcommand{\NS}{{\operatorname{NS}}}
\newcommand{\id}{{\operatorname{id}}}
\newcommand{\Zar}{{\text{\rm Zar}}} 
\newcommand{\Ord}{{\mathbf{Ord}}}
\newcommand{\FSimp}{{{\mathcal FS}}}
\newcommand{\inj}{{\text{\rm inj}}}  
\newcommand{\Sch}{{\operatorname{\mathbf{Sch}}}} 
\newcommand{\cosk}{{\operatorname{\rm cosk}}} 
\newcommand{\sk}{{\operatorname{\rm sk}}} 
\newcommand{\subv}{{\operatorname{\rm sub}}}
\newcommand{\bary}{{\operatorname{\rm bary}}}
\newcommand{\Comp}{{\mathbf{SC}}}
\newcommand{\IComp}{{\mathbf{sSC}}}
\newcommand{\Top}{{\mathbf{Top}}}
\newcommand{\holim}{\mathop{{\rm holim}}}
\newcommand{\Holim}{\mathop{{\it holim}}}
\newcommand{\op}{{\text{\rm op}}}
\newcommand{\<}{\mathopen<}
\renewcommand{\>}{\mathclose>}
\newcommand{\Sets}{{\mathbf{Sets}}}
\newcommand{\del}{\partial}
\newcommand{\fib}{{\operatorname{\rm fib}}}
\renewcommand{\max}{{\operatorname{\rm max}}}
\newcommand{\bad}{{\operatorname{\rm bad}}}
\newcommand{\Spt}{{\mathbf{Spt}}}
\newcommand{\Spc}{{\mathbf{Spc}}}
\newcommand{\Sm}{{\mathbf{Sm}}}
\newcommand{\cofib}{{\operatorname{\rm cofib}}}
\newcommand{\hocolim}{\mathop{{\rm hocolim}}}
\newcommand{\Glu}{{\mathbf{Glu}}}
\newcommand{\can}{{\operatorname{\rm can}}}
\newcommand{\Ho}{{\mathbf{Ho}}}
\newcommand{\GL}{{\operatorname{\rm GL}}}
\newcommand{\sq}{\square}
 \newcommand{\Ab}{{\mathbf{Ab}}}
\newcommand{\Tot}{{\operatorname{\rm Tot}}}
\newcommand{\loc}{{\operatorname{\rm s.l.}}}
\newcommand{\HZ}{{\operatorname{\sH \Z}}}
\newcommand{\Cyc}{{\operatorname{\rm Cyc}}}
\newcommand{\cyc}{{\operatorname{\rm cyc}}}
 \newcommand{\RCyc}{{\operatorname{\widetilde{\rm Cyc}}}}
 \newcommand{\Rcyc}{{\operatorname{\widetilde{\rm cyc}}}}
\newcommand{\Sym}{{\operatorname{\rm Sym}}}
\newcommand{\fin}{{\operatorname{\rm fin}}}
\newcommand{\SH}{{\operatorname{\sS\sH}}}
\newcommand{\Anna}{{$\mathcal{ANNA\ SARAH\ LEVINE}$}}
\newcommand{\Wedge}{{\Lambda}}
\newcommand{\eff}{{\operatorname{\rm eff}}}
\newcommand{\rcyc}{{\operatorname{\rm rev}}}
\newcommand{\DM}{{\operatorname{\mathcal{DM}}}}
\newcommand{\GW}{{\operatorname{\rm{GW}}}}
\newcommand{\sSets}{{\mathbf{Spc}}}
\newcommand{\Nis}{{\operatorname{Nis}}}
\newcommand{\et}{{\text{\'et}}}
\newcommand{\Cat}{{\mathbf{Cat}}}
\newcommand{\ds}{{/\kern-3pt/}}
\newcommand{\bD}{{\mathbf{D}}}
\newcommand{\bC}{{\mathbf{C}}}
\newcommand{\res}{{\operatorname{res}}}

\newcommand{\qf}{{\operatorname{q.fin.}}}

\newcommand{\fil}{\phi}
\newcommand{\barfil}{\sigma}

\newcommand{\Fac}{{\mathop{\rm{Fac}}}}
\newcommand{\Fun}{{\mathbf{Func}}}

\newcommand{\cp}{\coprod}
\newcommand{\ess}{\text{\rm{ess}}}
\newcommand{\dimrel}{\operatorname{\text{dim-rel}}}

 \newcommand{\gen}{{\text{gen}}}
\newcommand{\Gr}{{\text{\rm Gr}}}
\newcommand{\Ind}{{\operatorname{Ind}}}
\newcommand{\Supp}{{\operatorname{Supp}}}
\newcommand{\Cone}{{\operatorname{Cone}}}
\newcommand{\Tor}{{\operatorname{Tor}}}
\newcommand{\cok}{{\operatorname{coker}}}
\newcommand{\Br}{{\operatorname{Br}}}
\newcommand{\sm}{{\operatorname{sm}}}
\newcommand{\Proj}{{\operatorname{Proj}}}
\newcommand{\Pos}{{\operatorname{\mathbf{Pos}}}}
\newcommand{\colim}{\mathop{\text{colim}}}
\newcommand{\ps}{{\operatorname{\text{p.s.s.}}}}
\newcommand{\Pro}{{\mathbf{Pro}\hbox{\bf -}}}
\newcommand{\hocofib}{{\operatorname{\rm hocofib}}}
\newcommand{\hofib}{{\operatorname{\rm hofib}}}
\newcommand{\Th}{{\operatorname{\rm Th}}}
\newcommand{\Cube}{{\mathbf{Cube}}}
\newcommand{\Mod}{{\mathbf{Mod}}}
\newcommand{\alt}{{\operatorname{\rm alt}}}
\newcommand{\Alt}{{\operatorname{\rm Alt}}}
\newcommand{\sgn}{{\operatorname{\rm sgn}}}
\newcommand{\dga}{{\rm d.g.a.\!}}
\newcommand{\cdga}{{\rm c.d.g.a.\!}}
\newcommand{\sym}{{\operatorname{\rm sym}}}
\newcommand{\DTM}{{\operatorname{DTM}}}
\newcommand{\DCM}{{\operatorname{DCM}}}
\newcommand{\MTM}{{\operatorname{MTM}}}
\renewcommand{\sp}{{\widetilde{sp}}}
\newcommand{\NST}{{\mathbf{NST}}}
\newcommand{\PST}{{\mathbf{PST}}}
\newcommand{\Sus}{{\operatorname{Sus}}}
\newcommand{\tr}{{\operatorname{tr}}}
\newcommand{\cone}{{\operatorname{cone}}}
\newcommand{\na}{{\rm{naive}}}
\newcommand{\Cor}{{\operatorname{\mathbf{Cor}}}}
\newcommand{\MGL}{{\operatorname{MGL}}}

\addtocounter{section}{-1}
\title{Motivic tubular neighborhoods}

\author{Marc Levine}

\address{
Department of Mathematics\\
Northeastern University\\
Boston, MA 02115\\
USA}

\email{marc@neu.edu}

\keywords{}

\subjclass{Primary 14C25; Secondary 55P42, 18F20, 14F42}
\thanks{The author gratefully acknowledges the support of the Humboldt Foundation through the Wolfgang Paul Program, and support of the NSF via grants DMS 0140445 and  DMS-0457195}

\renewcommand{\abstractname}{Abstract}
\begin{abstract}  We construct  motivic versions of the classical tubular neighborhood and the punctured tubular neighborhood, and give applications to the construction of tangential base-points for mixed Tate motives, algebraic gluing of curves with boundary components, and limit motives.
\end{abstract}

\maketitle
\tableofcontents
\section{Introduction} Let $i:A\to B$ be a closed embedding of finite CW complexes. One useful fact is that $A$ admits a cofinal system of neighborhoods $T$ in $B$ with $A\to T$ a deformation retract. This is often used in the case where $B$ is a differentiable manifold, showing for example that $A$ has the homotopy type of the differentiable manifold $T$. This situation occurs in algebraic geometry, for instance in the case of a degeneration of smooth varieties $\sX\to C$ over the complex numbers.

To some extent, one has been able to mimic this construction in purely algebraic terms. The rigidity theorems of Gillet-Thomason \cite{GilletThomason}, extended by Gabber (details appearing a paper of Fujiwara \cite{Gabber}) indicated that, at least through the eyes of torsion \'etale sheaves, the topological tubular neighborhood can be replaced by the Hensel neighborhood. However, basic examples of non-torsion phenomena, even in the \'etale topology, show that the Hensel neighborhood cannot always be thought of as a tubular neighborhood, perhaps the simplest example being the sheaf $\G_m$.

Our object in this paper to to construct an algebraic version of the tubular neighborhood which has the basic properties of the topological construction, at least for a reasonably large class of cohomology theories. It turns out that a ``homotopy invariant" version of the Hensel neighborhood does the job, at least for theories which are themselves homotopy invariant. If one requires in addition that the given cohomology theory has a Mayer-Vietoris property for the Nisnevic topology, then one also has  an algebraic version of the punctured tubular neighborhood.

Our basic constructions are valid for a closed embedding of smooth varieties over a field. We extend these to the case of a (reduced) strict normal crossing subscheme by a Mayer-Vietoris procedure. Using an algebraic cosimplicial model of the universal cover, we can give an algebraic version of ``limit cohomology" for a semi-stable degeneration. We also have a monodromy sequence for theories with $\Q$-coefficients, and that satisfy an  additional restriction that the theory be ``alternating" (see Definition~\ref{Def:alt}). This allows us to give a definition of the limit motive of a semi-stable degeneration, and a motivic monodromy sequence with $\Q$-coefficients. We conclude with an application of our constructions to moduli of smooth curves and a construction of a specialization functor for category of mixied Tate motives, which in some cases yields a purely algebraic construction of tangential base-points.

We have left to another paper the task of checking the compatibilities of our constructions with others. Our punctured tubular neighborhood construction should be the same as some version of the functor $i^*Rj_*$ for the situation of a normal crossing scheme $i:D\to X$ with open complement $j:X\setminus D\to X$. Similarly, our limit cohomology construction should be a version of the sheaf of vanishing cycles, at least in the case of a semi-stable degeneration, and should be comparable with the constructions of Rappaport-Zink \cite{RappaportZink} as well as the limit mixed Hodge structure of Katz \cite{Katz} and Steenbrink \cite{Steenbrink}. Our specialization functor for Tate motives should be compatible with the Betti, \'etale and Hodge realizations; similarly, realization functors applied to our limit motive should yield for example the limit mixed Hodge structure. Additionally, Voevodsky's general setting of cross functors \cite{VoevCross} has now been shown (by Ayoub \cite{Ayoub} and R\"ondigs \cite{Roendigs}) to be available in the $\A^1$-stable homotopy category of $T$-spectra, and we expect our constructions are a weak (but also somewhat more concrete) version of some of the sheaf-theoretic constructions given by the general machinery. We have not attempted an investigation of this issue in this paper; however, we want to note that the construction given here for the limit cohomology in a semi-stable degeneration, and the monodromy sequence for limit cohomology, are taken directly from Ayoub's beautiful constructions of vanishing cycles and monodromy sequence in the setting of the $\A^1$-stable homotopy theory  of $T$-spectra.

My interest in this topic began as a result of several discussions on limit motives with Spencer Bloch and H\'el\`ene Esnault, whom I would like to thank for their encouragement and advice. I would also like to thank H\'el\`ene Esnault for clarifying the role of  the weight filtration leading to the exactness of Clemens-Schmidt monodromy sequence (see Remark~\ref{rem:MonodromySeq}). An earlier version of our constructions used an analytic (i.e. formal power series) neighborhood instead of the Hensel version now employed; I am grateful to Fabien Morel for suggesting this improvement. Finally, I want to thank Joseph Ayoub for explaining his construction of the vanishing cycles functor; his comments suggested to us  the use of the cosimplicial  path space in our construction of limit cohomology. In addition, Ayoub noticed a serious error in our first attempt at constructing the monodromy sequence; the method used in this version is following his suggestions and comments.

\section{Model structures and other preliminaries}
\subsection{Presheaves of simplicial sets}\label{sec:PreshSimp}  We recall some facts on the model structures in categories of simplicial sets, spectra, associated presheaf categories and certain localizations. For details, we refer the reader to \cite{Jardine1} and \cite{Hovey}.

We give $\sSets$ (simplicial sets) and $\sSets*$ (pointed simplicial sets) the standard model structures: cofibrations are (pointed) monomorphisms, weak equivalences are weak equivalences on the geometric realization, and fibrations are detemined by the right lifting property (RLP) with respect to trivial cofibrations; the fibrations are then exactly the Kan fibrations. We let  $|A|$ denote the geometric realization, and $[A,B]$  the homotopy classes of (pointed) maps $|A|\to |B|$.

We give  $\Spc(\sC)$  the model structure of functor categories described by Bousfield-Kan
 \cite{BousfieldKan}. That is, the cofibrations and weak equivalences are the pointwise ones, and the fibrations are determined by the RLP with respect to trivial cofibrations. We let $\sH\Spc(\sC)$   denote the associated homotopy category (see \cite{Hovey} for details).

\subsection{Presheaves of spectra} Let $\Spt$ denote the category of
spectra. To fix ideas, a spectrum will be a sequence of pointed simplicial sets
$E_0, E_1,\ldots$ together with maps of pointed simplicial sets $\epsilon_n:S^1\wedge E_n\to E_{n+1}$.  Maps of spectra are maps of the underlying simplicial sets which are compatible with the attaching maps $\epsilon_n$. The stable homotopy groups $\pi_n^s(E)$ are defined by
\[
\pi_n^s(E):=\lim_{m\to\infty} [S^{m+n},E_m].
\]

The category $\Spt$ has the following model structure: Cofibrations are maps $f:E\to F$ such that $E_0\to F_0$ is a cofibration, and for each $n\ge0$, the map
\[
E_{n+1}\coprod_{S^1\wedge E_n}S^1\wedge F_n\to F_{n+1}
\]
is a cofibration.  Weak equivalences are the stable weak equivalences, i.e., maps $f:E\to F$ which induce an isomorphism on $\pi_n^s$ for all $n$. Fibrations are characterized by having the RLP with respect to trivial cofibrations.

Let $\sC$ be a category. We say that a natural transformation $f:E\to E'$ of functors  $\sC^\op\to\Spt$ is a weak equivalence if $f(X):E(X)\to E'(X)$ is a stable weak equivalence for all $X$.  

We will assume that the category $\sC$ has an initial object $\0$ and admits finite coproducts over $\0$, denoted $X\amalg Y$. A functor $E:\sC^\op\to \Spt$ is called {\em additive} if for each $X, Y$ in $\sC$, the canonical map
\[
E(X\amalg Y)\to E(X)\oplus E(Y)
\]
is a weak equivalence. An additive functor  $E:\sC^\op\to \Spt$ is called a {\em presheaf of spectra} on $\sC$. This forms a full subcategory of the functor category.

We use the following model structure on the category of presheaves of spectra (see \cite{Jardine1}): Cofibrations and weak equivalences are given pointwise, and fibrations are characterized by having the RLP with respect to trivial cofibrations. We denote this model category by   $\Spt(\sC)$, and the associated homotopy category by $\sH\Spt(\sC)$.   We write $\SH$ for the homotopy category of $\Spt$.

Let $B$ be  a noetherian separated scheme of finite Krull dimension. We let $\Sm/B$ denote the category of smooth $B$-schemes of finite type over $B$. We often write $\Spt(B)$ and $\sH\Spt(B)$ for $\Spt(\Sm/B)$ and $\sH\Spt(\Sm/B)$.

For $Y\in\Sm/B$, a subscheme $U\subset Y$ of the form $Y\setminus \cup_\alpha F_\alpha$, with $\{F_\alpha\}$ a possibly infinite set of closed subsets of $Y$, is called {\em essentially smooth over } $B$; the category of essentially smooth $B$-schemes is denoted $\Sm^\ess$.

\subsection{Local model structure} 

If the category $\sC$ has a topology, there is often another model structure on $\Spc(\sC)$ or $\Spt(\sC)$ which takes this into account. We consider the case of the small Nisnevic site $X_\Nis$ on a scheme $X$ (assumed to be noetherian, separated and of finite Krull dimension), and on the big Nisnevic sites $\Sm/B_\Nis$ or $\Sch/B_\Nis$, as well as the Zariski versions $X_\Zar$, etc. We describe the Nisnevic version for spectra below; the definitions and results for the Zariski topology and for spaces are exactly parallel.

\begin{Def} A map $f:E\to F$ of presheaves of   spectra on $X_\Nis$ is a {\em local weak equivalence} if the induced map on the Nisnevic sheaf of stable homotopy groups $f_*:\pi_m^s(E)_\Nis\to \pi_m^s(F)_\Nis$ is an isomorphism of sheaves for all $m$. A map
$f:E\to F$ of presheaves of   spectra on $\Sm/B_\Nis$ or $\Sch/B_\Nis$ is a local weak equivalence if the restriction of $f$ to $X_\Nis$ is a local weak equivalence for all $X\in \Sm/B$ or $X\in\Sch/B$.
\end{Def}

The (Nisnevich) local model structure on the category of presheaves of spectra on $X_\Nis$ has cofibrations given pointwise, weak equivalences the local weak equivalences and  fibrations are characterized by having the RLP with respect to trivial cofibrations. We write $\Spt(X_\Nis)$ for  this model category,  and  $\sH\Spt(X_\Nis)$ for the associated homotopy category. The Nisnevic local model categories $\Spt_\Nis(\Sm/B)$ and $\Spt_\Nis(\Sch/B)$, with homotopy categories $\sH\Spt_\Nis(\Sm/B)$ and $\sH\Spt_\Nis(\Sch/B)$, are defined similarly. For details, we refer the reader to \cite{Jardine1}.

\subsection{$\A^1$-local structure} One can perform a Bousfield localization on  $\Spc(\Sm/B_\Nis)$ or $\Spt_\Nis(Sm/B)$ so that the maps $\Sigma^\infty X\times\A^1_+\to\Sigma^\infty X_+$ induced by the projections $X\times\A^1\to X$ become weak equivalences. We call the resulting model structure the {\em Nisnevic-local $\A^1$-model structure}, denoted $\Spc_{\A^1}(\Sm/B_\Nis)$ or $\Spt_{\A^1}(Sm/B_\Nis)$. One has Zariski-local and \'etale-local versions as well. We denote the homotopy categories for the Nisnevic version by $\sH_{\A^1}(B)$ (for spaces) and $\SH_{\A^1}(B)$ (for spectra). For the Zariski or \'etale versions, we indicate the topology in the notation.   For details, see \cite{MorelLec, Morel, MorelVoev}.

\subsection{Additional notation}   Given $W\in\Sm/S$, we have restriction functors 
\begin{align*}
&\Spc(S)\to\Spc(W_\Zar)\\
&\Spt(S)\to\Spt(W_\Zar);
\end{align*}
we write the restriction of some $E\in \Spc(S)$  to $\Spc(W_\Zar)$ as $E(W_\Zar)$. We use a similar notation for the restriction of $E$ to $\Spt(W_\Zar)$, or for restrictions to $W_\Nis$ or $W_\et$. More generally, if $p:Y\to W$ is a morphism in $\Sm/S$, we write $E(Y/W_\Zar)$ for the presheaf $U\mapsto E(Y\times_WU)$.

 For $Z\subset Y$ a closed subset, $Y\in\Sm/S$ and for $E\in \Spc(S)$ or $E\in\Spt(S)$, we write $E^Z(Y)$ for the homotopy fiber of the restriction map $E(Y)\to E(Y\setminus Z)$. We define the presheaf $E^{Z_\Zar}(Y)$ by setting, for $U\subset Z$ a Zariski open subscheme with closed complement $F$, 
\[
E^{Z_\Zar}(Y)(U):=E^U(Y\setminus F).
\]

A {\em co-presheaf} on a category $\sC$ with values in $\sA$ is just an $\sA$-valued preheaf on $\sC^\op$.

As usual, we let $\Delta^n$ denote the algebraic $n$-simplex
\[
\Delta^n:=\Spec \Z[t_0,\ldots,t_n]/\sum_it_i-1, 
\]
and $\Delta^*$ the cosimplicial scheme $n\mapsto \Delta^n$. For a scheme $X$, we have $\Delta^n_X:=X\times\Delta^n$ and the cosimplicial scheme $\Delta^*_X$.

Let $k$ be a field. For $E\in\Spc(k)$ or in $\Spt(k)$, we say that $E$ is {\em homotopy invariant} if for all $X\in\Sm/k$, the pull-ack map $E(X)\to E(X\times\A^1)$ is a   weak equivalence (resp., stable weak equivalence). We say that $E$ {\em satisfies Nisnevic excision} if $E$ transforms the standard Nisnevic squares (see 
\cite[Definition 1.3, pg. 96]{MorelVoev}) to homotopy cartesian  squares.

\section{Tubular neighborhoods for smooth pairs} 

Let $i:W\to X$ be a closed embedding in $\Sm/k$.  In this section, we construct the tubular neighborhood of $W$ in $X$ as a functor from $W_\Zar$ to cosimplicial pro-$k$-schemes, $\tau_\epsilon(W)$. Applying this functor to a space over $k$, $E\in\Spc(k)$ yields a presheaf of spaces $E(\tau_\epsilon(W))$ on $W_\Zar$, which is our main object of study.

\subsection{The cosimplicial pro-scheme $\tau_\epsilon(Z)$} For a closed embedding  $W\to T$ in $\Sm/k$, let $T^W_\Nis$ be the category of Nisnevic neighborhoods of $W$ in $T$, i.e., objects are \'etale maps $p:T'\to T$ of finite type, together with a section $s:W\to T'$ to $p$ over $W$. Morphisms are morphisms over $T$ which respect the sections. Note that $T^W_\Nis$ is a left-filtering essentially small category. 

Sending $(p:T'\to T, s:W\to T')$ to $T'\in \Sm/k$ defines the pro-object $\hat{T}^h_W$ of $\Sm/k$; the sections $s:W\to T'$ give rise to a map of the constant pro-scheme $W$ to $\hat{T}^h_W$, denoted
\[
\hat{i}_W:W\to \hat{T}^h_W. 
\]

Given a $k$-morphism $f:S\to T$, and closed embeddings $i_V:V\to S$, $i_W:W\to T$ such that $f\circ i_V$ factors through $i_W$ (by $\bar{f}:V\to W$), we have the pull-back functor
\begin{align*}
&f^*:T^W_\Nis\to S^V_\Nis,\\
&f^*(T'\to T, s:W\to T'):=(T'\times_TS, (s\circ \bar{f},i_V)).
\end{align*}
This gives us the map of pro-objects $f^h:\hat{S}^h_V\to \hat{T}^h_W$, so that sending $W\to T$ to $\hat{T}^h_W$ and $f$ to $f^h$ becomes a pseudo-functor. 

We let $f^h:\hat{S}^h_V\to \hat{T}^h_W$ denote the induced map on pro-schemes. If $f$ happens to be a Nisnevic neighborhood of $W\to X$ (so $\bar{f}:V\to W$ is an isomorphism) then $f^h:\hat{S}^h_V\to \hat{T}^h_W$ is clearly an isomorphism.

\begin{rem} The pseudo-functor $(W\to T)\mapsto \hat{T}^h_W$ can be rectified to an honest functor by first replacing  $T^W_\Nis$ with the cofinal subcategory $T^W_{\Nis,0}$ of neighborhoods $T'\to T$, $s:W\to T'$ such that each connected component of $T'$ has non-empty intersection with $s(W)$. One notes that  $T^W_{\Nis,0}$ has only identity automorphisms, so we replace $T^W_{\Nis,0}$  with a choice of a full subcategory $T^W_{\Nis,00}$  giving a set of representatives of the isomorphism classes in $T^W_{\Nis,0}$, We thus have the honest functor  $(W\to T)\mapsto T^W_{\Nis,00}$ which yields an equivalent pro-object $\hat{T}^h_W$. We will use the strictly functorial version from now on without comment.
\end{rem}

For a closed embedding $i:W\to X$ in $\Sm/k$, set $\hat{\Delta}^{n}_{X,W}:=(\widehat{\Delta^n_X})^h_{\Delta^n_W}$. The cosimplicial scheme 
\begin{align*}
&\Delta_X^*:\Ord\to\Sm/k\\
&[n]\mapsto \Delta^n_X
\end{align*}
thus gives rise to the cosimplicial pro-scheme 
\begin{align*}
&\hat{\Delta}_{X,W}^*:\Ord\to\Pro\Sm/k\\
&[n]\mapsto \hat{\Delta}^n_{X,W}
\end{align*}
The maps $\hat{i}_{\Delta^n_W}:\Delta^n_W\to (\widehat{\Delta^n_X})^h_{\Delta^n_W}$ give the closed embedding  of cosimplicial pro-schemes
\[
\hat{i}_W:\Delta_W^*\to \hat{\Delta}^n_{X,W}.
\]
Also, the canonical maps  $\pi_n:\hat{\Delta}^n_{X,W}\to \Delta^n_X$ define the map
\[
\pi_{X,W}:\hat{\Delta}^*_{X,W}\to\Delta^*_X.
\]

Let $(p:X'\to X, s:W\to X')$ be a Nisnevic neighborhood of $(W,X)$. The map 
\[
p:\hat{\Delta}^n_{X',W}\to \hat{\Delta}^n_{X,W}
\]
is an isomorphism respecting the embeddings $\hat{i}_W$. Thus, sending a Zariski open subscheme $U\subset W$ with complement $F\subset W\subset X$ to $\hat{\Delta}^n_{X\setminus F,U}$ defines a co-presheaf $\hat{\Delta}^n_{\hat{X},W_\Zar}$ on $W_\Zar$ with values in  pro-objects of $\Sm/k$; we write $\tau_\epsilon^{\hat{X}}(W)$ for the cosimplicial object
\[
n\mapsto \hat{\Delta}^n_{\hat{X},W_\Zar}. 
\]
We use the notation  $\hat{X}$ in the notation because the co-presheaf  $\hat{\Delta}^n_{\hat{X},W_\Zar}$ is independent of the choice of Nisnevic neighborhood $X$ of $W$, up to canonical isomorphism. 

Let $\Delta^*_{W_\Zar}$ denote the co-presheaf on $W_\Zar$ defined by $U\mapsto \Delta^*_U$. The closed embeddings $\hat{i}_U$ define the natural transformation 
\[
\hat{i}_W:\Delta^*_{W_\Zar}\to \tau_\epsilon^{\hat{X}}(W).
\]
The maps $\pi_{V,W\cap V}$ for $V\subset X$ a Zariski open subscheme define the map of functors
\[
\pi_{X,W}: \tau_\epsilon^{\hat{X}}(W)\to \Delta^*_{X_\Zar}
\]

\subsection{Evaluation on spaces}
Let $i:W\to T$ be a closed embedding in $\Sm/$. For $E\in\Spc(T)$, we  have the space $E(\hat{T}_h^W)$, defined by
\[
E(\hat{T}^h_W):=\colim_{(p:T'\to T, s:W\to T')\in T^W_\Nis}E(T').
\]
Given a Nisnevic neighborhood $(p:T'\to T, s:W\to T')$, we have the weak equivalence
\[
p^*:E(\hat{T}^h_W)\to E(\hat{T'}^h_{s(W)}).
\]
Thus, for each open subscheme $j:U\to W$, we may evaluate $E$ on the cosimplicial pro-scheme $\tau_\epsilon^{\hat{X}}(W)(U)$, giving us the presheaf of simplicial spectra $E(\tau_\epsilon^{\hat{X}}(W))$ on $W_\Zar$:
\[
E(\tau_\epsilon^{\hat{X}}(W))(U):=E(\tau_\epsilon^{\hat{X}}(W)(U)).
\]

Now suppose that $E$ is in $\Spc(k)$. 
The map $\hat{i}_W:\Delta^*_{W_\Zar}\to \tau_\epsilon^{\hat{X}}(W))$ gives us the  map of presheaves on $W_\Zar$
\[
i_W^*:E(\tau_\epsilon^{\hat{X}}(W))\to E(\Delta^*_{W_\Zar}).
\]
Similarly, the map $\pi_{X,W}$ gives the map of presheaves on $W_\Zar$
\[
\pi_{X,W}^*:E(\Delta^*_{X_\Zar})\to E(\tau_\epsilon^{\hat{X}}(W)).
\]

The main result of this section is
\begin{thm} \label{thm:main1} Let $E$ be in $\Spc(k)$. Then the map $i_W^*:E(\tau_\epsilon^{\hat{X}}(W))\to E(\Delta^*_{W_\Zar})$ is a weak equivalence for the Zariski-local model structure, i.e., for each point $w\in W$, the map $i_{W,w}^*$ on the stalks at $w$ is a weak equivalence of the associated total space.
\end{thm}

\subsection{Proof of Theorem~\ref{thm:main1}}
The proof relies on two lemmas. 

\begin{lem}\label{lem:A1Homotopy} Let $i:W\to X$ be a closed embedding in $\Sm/k$, giving the closed embedding $\A^1_W\to \A^1_X$. We have as well the maps $i_0, i_1:W\to \A^1_W$. Then for each $E\in\Spc(k)$, the maps
\[
i_0^*, i_1^*: E(\hat{\Delta}^*_{\A^1_X,\A^1_W})\to E(\hat{\Delta}^*_{X,W})
\]
are homotopic.
\end{lem}

\begin{proof} This is just an adaptation of the standard triangulation argument. For each order-preserving map $g=(g_1, g_2):[m]\to [1]\times[n]$, let
\[
T_g:\Delta^m\to\Delta^1\times \Delta^n,
\]
be the affine-linear extension of the map on the vertices 
\[
v_i\mapsto (v_{g_1(i)},v_{g_2(i)}). 
\]
$\id_X\times T_g$ induces the map
\[
\hat{T}_g:\hat{\Delta}^m_{X,W}\to (\widehat{\Delta^1\times \Delta^n_X})^h_{\Delta^1\times \Delta^n_W}
\]
We note that the isomorphism $(t_0,t_1)\mapsto t_0$ of $(\Delta^1,v_1,v_0)$ with $(\A^1,0,1)$ induces an isomorphism of cosimplicial schemes
\[
\hat{\Delta}^*_{\A^1_X,\A^1_W}\cong  (\widehat{\Delta^1\times \Delta^*_X})^h_{\Delta^1\times \Delta^*_W}.
\]

The maps 
\[
\hat{T}^*_g:E(\hat{\Delta}^n_{\A^1_X,\A^1_W})\to E(\hat{\Delta}^m_{X,W})
\]
induce a simplicial homotopy $T$ between $i_0^*$ and $i_1^*$: indeed, the maps $T_g$ satisfy the identities necessary to define a map of cosimplicial schemes
\[
T:\Delta^*\to (\Delta^1\times \Delta^*)^{\Delta[1]},
\]
with $\delta_0^*\circ T=i_0$, $ \delta_1^*\circ T=i_1$. Applying the functor ${}^h$, we see that the maps $\hat{T}_g$ define the map of cosimplicial schemes
\[
\hat{T}:\hat{\Delta}^*_{X,W}\to (\hat{\Delta}^*_{\A^1_X,\A^1_W})^{\Delta[1]},
\]
with $\hat{T}\circ \delta_0=\hat{i}_0$, $\hat{T}\circ\delta_1=\hat{i}_1$; we then apply $E$.
\end{proof}

\begin{lem} \label{lem:TrivEquiv}Take $W\in\Sm_k$. Let $X=\A^n_W$ and let $i:W\to X$ be the 0-section. Then $i_W^*:E(\hat{\Delta}^*_{X,W})\to E(\Delta^*_W)$ is a homotopy equivalence.
\end{lem}

\begin{proof} Let $p:X\to W$ be the projection, giving the map 
\[
\hat{p}: \hat{\Delta}^*_{X,W}\to
\hat{\Delta}^*_{W,W}=\Delta^*_W
\]
 and  $\hat{p}^*:E(\Delta^*_W)\to E(\hat{\Delta}^n_{X,W})$. Clearly $\hat{i}_W^*\circ \hat{p}^*=\id$, so it suffices to show that $\hat{p}^*\circ \hat{i}_W^*$ is homotopic to the identity.

For this, we use the  multiplication map  $\mu:\A^1\times\A^n\to \A^n$,
\[
\mu(t;x_1,\ldots, x_n):=(tx_1,\ldots, tx_n).
\]
The map $\mu\times\id_{\Delta^*}$ induces the map
\[
\hat{\mu}:(\widehat{\A^1\times\A^n_W\times\Delta^*})^h_{\A^1\times 0_W\times\Delta^*}\to
(\widehat{\A^n_W\times\Delta^*})^h_{0_W\times\Delta^*}
\]
with $\hat{\mu}\circ \hat{i}_0=\hat{i}_W\circ \hat{p}$ and 
 $\hat{\mu}\circ\hat{i}_1 =\id$. Since $\hat{i}_0^*$ and $\hat{i}_1^*$ are homotopic by Lemma~\ref{lem:A1Homotopy} , the proof is complete.
 \end{proof}

To complete the proof of Theorem~\ref{thm:main1}, take a point $w\in W$. Then replacing $X$ with a Zariski open neighborhood of $w$, we may assume there is a Nisnevic neigborhood $X'\to X$, $s:W\to X'$ of $W$ in $X$ such that $W\to X'$ is in turn a Nisnevic neighborhood of the zero-section $W\to \A^n_W$, $n=\codim_XW$. Since $E(\hat{\Delta}^n_{X,W})$ is thus weakly equivalent  to $E(\hat{\Delta}^n_{\A^n_W,0_W})$, the result follows from Lemma~\ref{lem:TrivEquiv}.

\begin{cor}\label{cor:Retract} Suppose that $E\in\Spc(k)$ is fibrant for the Zariski-local, $\A^1$ model structure. Then for $i:W\to X$ a closed embedding,  there is a natural isomorphism in $\sH\Spc(W_\Zar)$
\[
E(\tau_\epsilon^{\hat{X}}(W))\cong E(\tau_\epsilon^{\hat{N_i}}(0_W))
\]
Here $N_i$ is the normal bundle of the embedding $i$, and $0_W\subset N_i$ is the 0-section.
\end{cor}
\begin{proof} This follows directly from Theorem~\ref{thm:main1}: Since $E$ is fibrant for the  Zariski-local, $\A^1$ model structure, $E$ is homotopy invariant. In consequence, for each $T\in\Sm/k$, the canonical map
\[
E(T)\to E(\Delta^*_T)
\]
is a weak equivalence. The desired isomorphism in $\sH\Spc(W_\Zar)$ is constructed by
composing the isomorphisms
\[
E(\tau_\epsilon^{\hat{X}}(W))\to E(\Delta^*_W)\leftarrow E(W)=E(0_W)\to E(\Delta^*_{0_W})\leftarrow
E(\tau_\epsilon^{\hat{N_i}}(0_W)).
\]
\end{proof}

\section{Punctured tubular neighborhoods}

The real interest is not in the tubular neighborhood $\tau_\epsilon^{\hat{X}}(W)$, but in the punctured tubular neighborhood $\tau_\epsilon^{\hat{X}}(W)^0$. In this section, we define this object and discuss its basic properties.

\subsection{Definition of the punctured neighborhood} Let $i:W\to X$ be a closed embedding in $\Sm/k$. We have the closed embedding of cosimplicial pro-schemes
\[
\hat{i}:\Delta^*_W\to \hat{\Delta}^*_{X,W};
\]
we let
\[
\hat{j}: \hat{\Delta}^*_{X\setminus W}\to  \hat{\Delta}^*_{X,W}
\]
be the open complement $\hat{\Delta}^n_{X\setminus W}:=\hat{\Delta}^n_{X,W}\setminus
\Delta^n_W$. Extending this construction to all open subschemes of $X$, we have the co-presheaf on $W_\Zar$, 
\[
U=W\setminus F \mapsto \hat{\Delta}^*_{(X\setminus F)\setminus U},
\]
which we denote by $\tau_\epsilon^{\hat{X}}(W)^0$. The maps 
\[
\hat{j}_U: \hat{\Delta}^*_{(X\setminus F)\setminus U}\to  \hat{\Delta}^*_{(X\setminus F), U}
\]
define the map $\hat{j}:\tau_\epsilon^{\hat{X}}(W)^0\to \tau_\epsilon^{\hat{X}}(W)$. The maps $\hat{\Delta}^*_{U\setminus W\cap U}\to \Delta^*_{U\setminus W\cap U}$ give us the map
\[
\pi:\tau_\epsilon^{\hat{X}}(W)^0\to \Delta^*_{(X\setminus W)_\Zar}.
\]

To give a really useful result on the presheaf $E(\tau_\epsilon^{\hat{X}}(W)^0)$, we will need to impose additional conditions on $E$. These are
\begin{enumerate}
\item $E$ is homotopy invariant
\item $E$ satisfies Nisnevic excision
\end{enumerate}

One important consequence of these properties is the purity theorem of Morel-Voevodsky:
\begin{thm}[Purity \hbox{\cite{MorelVoev}}] Suppose $E\in\Spt(k)$ is homotopy invariant and satisfies Nisnevic excision.  Then there is an isomorphism in $\sH\Spt(W_\Zar)$
\[
E^{W_\Zar}(X)\to E^{0_{W \Zar}}(N_i)
\]
\end{thm}

Let $E(X/W_\Zar)$ be the presheaf on $W_\Zar$
\[
W\setminus F\mapsto E(X\setminus F)
\]
and $E(X\setminus W)$ the constant preheaf.

Let $\res:E(X/W_\Zar)\to E(\tau_\epsilon^{\hat{X}}(W))$, $\res^0:E(X\setminus W)\to E(\tau_\epsilon^{\hat{X}}(W)^0)$ be the pull-back by the natural maps $\tau_\epsilon^{\hat{X}}(W)(W\setminus F)\to X\setminus F$, 
$\tau_\epsilon^{\hat{X}}(W)^0\to X\setminus W$. Let $E^{\Delta^*_W}(\tau_\epsilon^{\hat{X}}(W))\in\Spt(W_\Zar)$ be the homotopy fiber of
the restriction map 
\[
\hat{j}:E(\tau_\epsilon^{\hat{X}}(W))\to E(\tau_\epsilon^{\hat{X}}(W)^0).
\]
The commutative diagram in $\Spt(W_\Zar)$ 
\[
\xymatrix{
E(X/W_\Zar)\ar[r]^{j^*}\ar[d]_\res& E(X\setminus W)\ar[d]^{\res^0}\\
E(\tau_\epsilon^{\hat{X}}(W))\ar[r]_{\hat{j}^*}& E(\tau_\epsilon^{\hat{X}}(W)^0)
}
\]

induces the map of distinguished triangles
\[
\xymatrix{
E^{W_\Zar}(X)\ar[r]\ar[d]_\psi&E(X/W_\Zar)\ar[r]^{j^*}\ar[d]_\res& E(X\setminus W)\ar[d]^{\res^0}\\
E^{\Delta^*_W}(\tau_\epsilon^{\hat{X}}(W))\ar[r]&E(\tau_\epsilon^{\hat{X}}(W))\ar[r]_{j^*}& E(\tau_\epsilon^{\hat{X}}(W)^0)
}
\]

We can now state the main result for $E(\tau_\epsilon^{\hat{X}}(W)^0)$. 

\begin{thm}\label{thm:main2} Suppose that $E\in\Spt(k)$ is homotopy invariant and satisfies Nisnevic excision. Let $i:W\to X$ be a closed embedding in $\Sm/k$. Then the map $\psi$ is a Zariski local weak equivalence.
\end{thm}

\begin{proof}  
Let $i_{\Delta^*}:\Delta^*_W\to \Delta^*_X$ be the embedding $\id\times i$.  For $U=W\setminus F\subset W$, $\tau_\epsilon^{\hat{X}}(W)^0(U)$ is the cosimplicial scheme with $n$-cosimplices 
\[
\tau_\epsilon^{\hat{X}}(W)^0(U)^n= \hat{\Delta}^n_{X\setminus F,U}\setminus\Delta^n_U
\]
so we have the natural isomorphism 
\[
E^{\Delta^*_W}(\tau_\epsilon^{\hat{X}}(W))\cong E^{\Delta^*_{W_\Zar}}(\Delta^*_X/W_\Zar),
\]
where $E^{\Delta^*_{W_\Zar}}(\Delta^*_X/W_\Zar)(W\setminus F)$ is the total spectrum of the   simplicial spectrum
\[
n\mapsto E^{\Delta^n_{W\setminus F}}(\Delta^n_{X\setminus F}).
\]
The homotopy invariance of $E$ implies that the pull-back
\[
E^{W\setminus F}(X\setminus F)\to E^{\Delta^n_{W\setminus F}}(\Delta^n_{X\setminus F})
\]
is a weak equivalence for all $n$, so we have the weak equivalence
\[
E^{W_\Zar}(X)\to E^{\Delta^*_{W_\Zar}}(\Delta^*_X/W_\Zar).
\]
It follows from the construction that this is the map $\psi$, completing the proof.
\end{proof}

\begin{cor}\label{cor:DistTriang1} There is a distinguished triangle in 
$\sH\Spt(W_\Zar)$
\[
E^{W_\Zar}(X)\to E(W_\Zar)\to E(\tau_\epsilon^{\hat{X}}(W)^0)
\]
\end{cor}

\begin{proof} By Theorem~\ref{thm:main1}, the map $\hat{i}^*:E( \tau_\epsilon^{\hat{X}}(W))\to E(\Delta^*_{W_\Zar})$ is 
a weak equivalence; using homotopy invariance again, the map
\[
E(W_\Zar)\to E(\Delta^*_{W_\Zar})
\]
is a weak equivalence. Combining this with Theorem~\ref{thm:main2} yields the result.
\end{proof}

For homotopy invariant $E\in\Spt(k)$, we let 
\[
\phi_E:E(\tau_\epsilon^{\hat{N_i}}(0_W))\to  E(\tau_\epsilon^{\hat{X}}(W)).
\]
be the isomorphism in $\sH\Spt(W_\Zar)$ defined as the composition
\[
E(\tau_\epsilon^{\hat{N_i}}(0_W))\cong E(0_W)=E(W)\cong  E(\tau_\epsilon^{\hat{X}}(W)),
\]
where the isomorphisms in the line above are given by Theorem~\ref{thm:main1}.

\begin{cor}\label{cor:DistTriang2} Suppose that $E\in\Spt(k)$ is homotopy invariant and satisfies Nisnevic excision. Let $i:W\to X$ be a closed embedding in $\Sm/k$ and let $N_i^0=N_i\setminus 0_W$. 

\medskip
\noindent
(1) The restriction maps
\begin{align*}
&\res:E(N_i/W_\Zar)\to E(\tau_\epsilon^{\hat{N_i}}(0_W))\\
&\res^0:E(N_i^0/W_\Zar)\to E(\tau_\epsilon^{\hat{N_i}}(0_W)^0)
\end{align*}
are weak equivalences in $\Spt(W_\Zar)$.

\medskip
\noindent
(2)  There is a canonical isomorphism in $\sH\Spt(W_\Zar)$
\[
\phi^0_E:E(\tau_\epsilon^{\hat{N_i}}(0_W)^0)\to E(\tau_\epsilon^{\hat{X}}(W)^0)
\]

\noindent
(3)  Consider the diagram (in $\sH\Spt(W_\Zar)$)
\[
\xymatrix{
E^{0_{W\Zar}}(N_i)\ar[r]\ar@{=}[d]i& E(N_i/W_\Zar)\ar[d]^{\res_E}\ar[r] &E(N^0_i/W_\Zar)\ar[d]^{\res_E^0}\\
E^{0_{W\Zar}}(N_i)\ar[r]\ar[d]_\pi&E(\tau_\epsilon^{\hat{N_i}}(0_W))\ar[r]^{\hat{j_N}^*}\ar[d]^-{\phi_E}
&E(\tau_\epsilon^{\hat{N_i}}(0_W)^0)\ar[d]^-{\phi_E^0}\\
E^{W_\Zar}(X)\ar[r]& E(\tau_\epsilon^{\hat{X}}(W))\ar[r]_{\hat{j}^*}&E(\tau_\epsilon^{\hat{X}}(W)^0)\\
E^{W_\Zar}(X)\ar[r]\ar@{=}[u]& E(X/W_\Zar)\ar[u]_{\res_E}\ar[r]_{j^*}& E(X\setminus W)\ar[u]_{\res_E^0}
}
\]
The first and last rows are the homotopy fiber sequences defining the presheaves  $E^{0_{W\Zar}}(N_i)$ and $E^{W_\Zar}(X)$, respectively, the second row and third rows are the distinguished triangles of Theorem~\ref{thm:main2},   and $\pi$ is the Morel-Voevodsky purity isomorphism. Then this diagram commutes and each triple of vertical maps defines a map of distinguished triangles.
\end{cor}

\begin{proof} It follows directly from the weak equivalence (in Theorem~\ref{thm:main2}) of the homotopy fiber of 
\[
\hat{j}^*:E(\tau_\epsilon^{\hat{X}}(W))\to E(\tau_\epsilon^{\hat{X}}(W)^0)
\]
with $E^{W_\Zar}(X)$ that the triple $(\id, \res_E,\res_E^0)$ defines a map of distinguished triangles. The same holds for the map of the first row to the second row; we now verify that this latter map is an isomorphism of distinguished triangles.

For this, let $s:W\to N_i$ be the zero-section. We have the isomorphism $i_W^*:E(\tau_\epsilon^{\hat{N_i}}(0_W))\to E(W_\Zar)$  defined as the diagram of weak equivalences
\[
E(\tau_\epsilon^{\hat{N_i}}(0_W))\xrightarrow{i_{\Delta^*_{W_\Zar}}^*}
E(\Delta^*_{W_\Zar})\xleftarrow{\iota_{0*}}E(W_\Zar).
\]
From this, it is easy to check that the diagram
\[
\xymatrix{
E(N_i/W_\Zar)\ar[r]^{\res_E}\ar[dr]_{s^*}&E(\tau_\epsilon^{\hat{N_i}}(0_W))\ar[d]^{i_W^*}\\
&E(W_\Zar)
}
\]
commutes in  $\sH\Spt(W_\Zar)$. As $E$ is homotopy invariant, $s^*$ is an isomorphism, hence $\res_E$ is an isomorphism as well. This completes the proof of (1).

The proof of (2) and (3)  uses the standard deformation diagram. Let $\bar\mu:\bar Y\to X\times\A^1$ be the blow-up of $X\times\A^1$ along $W\time0$, let $\bar\mu^{-1}[X\times0]$ denote the proper transform, and let $\mu:Y\to X\times\A^1$ be the open subscheme $\bar Y\setminus\bar\mu^{-1}[X\times0]$. Let $p:Y\to \A^1$ be $p_2\circ \mu$. Then $p^{-1}(0)=N_i$, $p^{-1}(1)=X\times1=X$, and $Y$ contains the proper transform $\bar\mu^{-1}[W\times\A^1]$, which is mapped isomorphically by $\mu$ to $W\times\A^1\subset X\times\A^1$. We let $\tilde{i}i:W\times\A^1\to Y$ be the resulting closed embedding. The restriction of $\tilde{i}$ to $W\times0$ is the zero-section $s:W\to N_i$ and the restriction of $\tilde{i}$ to $W\times1$ is $i:W\to X$. The resulting diagram is
\begin{equation}\label{eqn:DefDiag}
\xymatrix{
W\ar[r]^-{i_0}\ar[d]_s&W\times\A^1\ar[d]_{\tilde i}&W\ar[l]_-{i_1}\ar[d]^i\\
N_i\ar[r]^-{i_0}\ar[d]_{p_0}&Y\ar[d]_p&X\ar[l]_-{i_1}\ar[d]^{p_1}\\
0\ar[r]_{i_0}&\A^1&1\ar[l]^{i_1}
}
\end{equation}

Together with Theorem~\ref{thm:main2}, diagram \eqref{eqn:DefDiag} gives us two maps of distinguished triangles:
\begin{multline*}
[E^{0_{W\times\A^1\Zar}}(Y)\to E( \tau_\epsilon^{\hat{Y}}(W\times\A^1))\to E(\tau_\epsilon^{\hat{Y}}(W\times\A^1)^0)]\\
\xrightarrow{i_1^*}\\
[E^{W_\Zar}(X)\to E( \tau_\epsilon^{\hat{X}}(W))\to E(\tau_\epsilon^{\hat{X}}(W)^0)]
\end{multline*}
and

\begin{multline*}
[E^{0_{W\times\A^1\Zar}}(Y)\to E(\tau_\epsilon^{\hat{Y}}(W\times\A^1))\to E(\tau_\epsilon^{\hat{Y}}(W\times\A^1)^0)]\\
\xrightarrow{i_0^*}\\
[E^{0_{W\Zar}}(N_i)\to E(\tau_\epsilon^{\hat{N_i}}(0_W))\to E(\tau_\epsilon^{\hat{N_i}}(0_W)^0)]
\end{multline*}

As above, we have the commutative diagram
\[
\xymatrix{
E( \tau_\epsilon^{\hat{Y}}(W\times\A^1))\ar[r]^{i_0^*}\ar[d]_{i_{W\times\A^1}^*}&
E( \tau_\epsilon^{\hat{N_i}}(0_W))\ar[d]_{i_W^*}\\
E(W\times\A^1)\ar[r]_{i_0^*}&E(W).
}
\]
As $E$ is homotopy invariant, the maps $i_W^*$, $i_{W\times\A^1}^*$ and $i_0^*:E(W\times\A^1)\to E(W)$ are isomorphisms, hence 
\[
i_0^*:E( \tau_\epsilon^{\hat{Y}}(W\times\A^1))\to E( \tau_\epsilon^{\hat{N_i}}(0_W))
\]
 is an isomorphism. Similarly, 
 \[
i_1^*:E( \tau_\epsilon^{\hat{Y}}(W\times\A^1))\to E( \tau_\epsilon^{\hat{X}}(W))
\]
 is an isomorphism. The proof of the Morel-Voevodsky purity theorem \cite[Theorem 2.23]{MorelVoev} shows that
 \begin{align*}
 &i_0^*:E^{0_{W\times\A^1\Zar}}(Y)\to E^{0_W}(N_i)\\
 &i_1^*:E^{0_{W\times\A^1\Zar}}(Y)\to E^{W_\Zar}(X)
 \end{align*}
 are weak equivalences; the isomorphism $\pi$ is by definition $i_1^*\circ(i_0^*)^{-1}$. Thus, both $i_0^*$ and $i_1^*$ define isomorphisms of distinguished triangles, and 
 \[
 i_1^*\circ(i_0^*)^{-1}:E( \tau_\epsilon^{\hat{N_i}}(0_W))\to
 E( \tau_\epsilon^{\hat{X}}(W))
 \]
 is the map $\phi_E$. Defining $\phi_E^0$ to be the isomorphism
 \[
 i_1^*\circ(i_0^*)^{-1}:E( \tau_\epsilon^{\hat{N_i}}(0_W)^0)\to
 E( \tau_\epsilon^{\hat{X}}(W)^0)
 \]
 proves both (2) and (3).
 \end{proof}
 
 \begin{rems}\label{rems:TubFunct} (1) It follows from the construction of $\phi_E$ and $\phi_E^0$ that both these maps are natural in $E$.\\
 \\
 (2) The maps  $\phi_E^0$ are functorial in the embedding $i:W\to X$ in the following sense: Suppose we have closed embeddings $i_1:W\to X$, $i_2:X\to Y$. Fix $E$ and let  $\phi_{jE}^0$ be the map corresponding to the embeddings $i_j$, $j=1,2$. Similarly, $\phi_{E}$, $\phi_{E}^0$ be the maps corresponding to the embedding $i:=i_2\circ i_1$. We have the evident maps
 \begin{align*}
& \iota:\tau_\epsilon^{\hat{X}}(W)^0\to \tau_\epsilon^{\hat{Y}}(W)^0
& \eta:\tau_\epsilon^{\hat{N_{i_1}}}(W)^0\to \tau_\epsilon^{\hat{N_{i_2}}}(W)^0
\end{align*}
Then the diagram
\[
\xymatrix{
E(\tau_\epsilon^{\hat{N_{i_2}}}(W)^0)\ar[r]^{\phi_{2E}}\ar[d]_{\eta^*}&E(\tau_\epsilon^{\hat{Y}}(W)^0)\ar[d]^{\iota^*}\\
E(\tau_\epsilon^{\hat{N_{i_1}}}(W)^0)\ar[r]_{\phi_{1E}}&E(\tau_\epsilon^{\hat{X}}(W)^0)
}
\]
commutes. This follows by considering the ``double deformation diagram" associated to the two embeddings $i_1, i_2$, as in the proof of the functoriality of the Gysin morphism in \cite[Chap. III, Proposition 2.2.1]{LevineMixMot}.
 \end{rems}
 
\subsection{The exponential map}\label{sec:exp} Let $i:W\to X$ be a closed embedding in $\Sm/k$ with normal bundle $N_i\to W$. We have the map 
\[
\exp:N_i\to X
\]
in $\Spc(k)$, defined as the composition $N_i\to W\to X$. We also have the Morel-Voevodsky purity isomorphism
\[
\pi:\Th(N_i)\to X/(X\setminus W)
\]
in $\sH(k)$, giving the commutative diagram in $\sH(k)$
\begin{equation}\label{eqn:PreExp}
\xymatrix{
N_i\ar[r] \ar[d]_{\exp}&\Th(N_i)\ar[d]^\pi\\
X\ar[r]&X/(X\setminus W)
}
\end{equation}

Since we have the homotopy cofiber sequences:
\begin{align*}
&N_i\setminus 0_W\to N_i\to \Th(N_i)\to \Sigma (N_i\setminus 0_W)_+\\
&X\setminus W\to X\to X/(X\setminus W)\to \Sigma(X\setminus W)_+
\end{align*}
the diagram \eqref{eqn:PreExp} induces a map
\[
\Sigma (N_i\setminus 0_W)_+\to \Sigma(X\setminus W)_+
\]
in $\sH(k)$, however, this map is not canonical.

In this section we define a canonical map
\[
\exp^0:\Sigma^\infty (N_i\setminus 0_W)_+\to \Sigma^\infty(X\setminus W)_+
\]
in $\SH_{\A^1}(k)$ which yields the map of distinguished triangles in $\SH_{\A^1}(k)$:
\[
\xymatrix{
\Sigma^\infty (N_i\setminus 0_W)_+\ar[d]_{\exp^0}\ar[r]&\Sigma^\infty N_{i+}\ar[r] \ar[d]_{\exp}&\Sigma^\infty\Th(N_i)\ar[d]^\pi\\
\Sigma^\infty(X\setminus W)_+\ar[r]&\Sigma^\infty X_+\ar[r]&\Sigma^\infty X/(X\setminus W)
}
\]

To define $\exp^0$, we apply Corollary~\ref{cor:DistTriang2}  with $E$ a fibrant model of $\Sigma^\infty (X\setminus W)_+$. Denote the composition
\begin{multline*}
E(X\setminus W)\xrightarrow{\res^0_E} E(\tau_\epsilon^{\hat{X}}(W)^0)(W)\\
\xrightarrow{(\phi^0_E)^{-1}}E(\tau_\epsilon^{\hat{N_i}}(0_W)^0)(W)
\xrightarrow{(\res^0_E)^{-1}}E(N_i^0)
\end{multline*}
by $\exp^{0*}_E$. Since $E$ is fibrant, we have canonical isomorphisms
\begin{align*}
\pi_0E(N_i^0)&\cong\Hom_{\SH_{\A^1}(k)}(\Sigma^\infty N_{i+}^0, E)\\
&\cong
\Hom_{\SH_{\A^1}(k)}(\Sigma^\infty N_{i+}^0, \Sigma^\infty (X\setminus W)_+)\\
\\
\pi_0E(X\setminus W)&\cong\Hom_{\SH_{\A^1}(k)}(\Sigma^\infty (X\setminus W)_+, E)\\
&\cong \Hom_{\SH_{\A^1}(k)}(\Sigma^\infty (X\setminus W)_+, \Sigma^\infty (X\setminus W)_+)
\end{align*}
so $\exp^{0*}_E$ induces the map
\begin{multline*}
\Hom_{\SH_{\A^1}(k)}(\Sigma^\infty (X\setminus W)_+, \Sigma^\infty (X\setminus W)_+)\\
\xrightarrow{\exp^{0*}_E}
\Hom_{\SH_{\A^1}(k)}(\Sigma^\infty N_{i+}^0, \Sigma^\infty (X\setminus W)_+).
\end{multline*}
We set
\[
\exp^0:=\exp^{0*}_E(\id).
\]

To finish the construction, we show
\begin{prop}\label{prop:exp} The diagram, with rows the evident homotopy cofiber sequences,
\[
\xymatrixcolsep{15pt}
\xymatrix{
\Sigma^\infty (N_i\setminus 0_W)_+\ar[d]_{\exp^0}\ar[r]&\Sigma^\infty N_{i+}\ar[r] \ar[d]_{\exp}&\Sigma^\infty\Th(N_i)\ar[d]^\pi\ar[r]^-\partial&\Sigma\Sigma^\infty (N_i\setminus 0_W)_+\ar[d]_{\Sigma\exp^0}\\
\Sigma^\infty(X\setminus W)_+\ar[r]&\Sigma^\infty X_+\ar[r]&\Sigma^\infty X/(X\setminus W)\ar[r]_-\partial
&\Sigma\Sigma^\infty(X\setminus W)_+
}
\]
commutes in $\SH_{\A^1}(k)$.
\end{prop}

\begin{proof} It suffices to show that, for all fibrant $E\in\Spt(k)$, the diagram formed by applying $\Hom_{\SH_{\A^1}(k)}(-,E)$ to our diagram commutes. This latter diagram is the same as applying $\pi_0$ to the diagram
\begin{equation}\label{eqn:ComDiag}
\xymatrix{
E(N_i\setminus 0_W)&E(N_i)\ar[l]&E^{0_W}(N_i)\ar[l]&\Omega E(N_i\setminus 0_W))\ar[l]_-\partial\\
E(X\setminus W)\ar[u]_{\exp^{0*}}&E(X) \ar[u]_{\exp^*}\ar[l]&E^W(X)\ar[u]^{\pi^*}\ar[l]
&\Omega E(X\setminus W)\ar[l]^-\partial\ar[u]_{\Omega\exp^{0*}}
}
\end{equation}
where the rows are the evident homotopy fiber sequences.  It follows by the definition of $\exp^0$ and $\exp$ that this diagram is just the ``outside" of the diagram in Corollary~\ref{cor:DistTriang2}(3), extended to make the distinguished triangles explicit. Thus the diagram \eqref{eqn:ComDiag} commutes, which finishes the proof.
\end{proof}

\begin{rem} \label{rem:expFunct} The exponential maps $\exp$ and $\exp^0$ are functorial with respect to composition of embeddings, in the sense of the functoriality of the maps $\phi_E$, $\phi_E^0$ described in Remark~\ref{rems:TubFunct}. This follows from the functoriality of the maps $\phi_E$, $\phi_E^0$.
\end{rem}

\section{Neighborhoods of normal crossing schemes} We extend our results to the case of a strict normal crossing divisor $W\subset X$ by using a Mayer-Vietoris construction.

\subsection{Normal crossing schemes} Let $D$ be a reduced effective Cartier divisor on a smooth $k$-scheme $X$ with irreducible components $D_1$,$\ldots$, $D_m$. For each $I\subset \{1,\ldots, m\}$, we set
\[
D_I:=\cap_{i\in I}D_i
\]
We let $|D|$ denote the support of $D$, $i:|D|\to X$ the inclusion. For each $I\neq\0$, we let $\iota_I:D_I\to|D|$, $i_I:D_I\to X$ be the inclusions; for $I\subset J\subset \{1,\ldots, m\}$ we have as well the inclusion $i_{I,J}:D_J\to D_I$. Let $D^{(n)}:=\cup_{|I|=n}D_I$, $\tilde{D}^{(n)}:=\amalg_{|I|=n}D_I$, $i_{(n)}:D^{(n)}\to X$, $\iota_{(n)}:D^{(n)}\to|D|$ the inclusions and $\pi_{(n)}:\tilde{D}^{(n)}\to|D|$, 
$p_{(n)}:\tilde{D}^{(n)}\to X$ the evident maps.

Recall that $D$ is a {\em strict normal crossing divisor} if each for each $I$,  $D_I$ is smooth over $k$ and $\codim_XD_I=|I|$.

We extend this notion a bit: We call a closed subscheme $D\subset X$ a {\em strict normal crossing subscheme} if $X$ is in $\Sm/k$ and there is a smooth locally closed subscheme $Y\subset X$ containing $D$ such that $D$ is a strict normal crossing divisor on $Y$

\subsection{The tubular neighborhood}

Let $D\subset X$ be a strict normal crossing subscheme with irreducible components $D_1,\ldots, D_m$. For each $I\subset\{1,\ldots, m\}$, $I\neq \0$, we have the tubular neighborhood co-presheaf $\tau_\epsilon^{\hat{X}}(D_I)$ on $D_I$; taking the disjoint union gives the co-presheaf  $\tau_\epsilon^{\hat{X}}(\tilde{D}^{(n)})$ on $D^{(n)}$. The various inclusions $i_{I,J}$ give us the maps of co-presheaves
\[
\hat{i}_{I,J}:i_{I,J*}(\tau_\epsilon^{\hat{X}}(D_J))\to \tau_\epsilon^{\hat{X}}(D_I);
\]
pushing forward by the maps $\iota_I$ yields the diagram of co-presheaves on $|D|$
\begin{equation}\label{eqn:CoMVDiag1}
I\mapsto \iota_{I*}(\tau_\epsilon^{\hat{X}}(D_I))
\end{equation}
indexed by the non-empty $I\subset\{1,\ldots, m\}$. We have as well the diagram of identity co-presheaves
\begin{equation}\label{eqn:CoMVDiag2}
I\mapsto \iota_{I*}(D_I)
\end{equation}
and the natural transformation $\hat{\iota}_{|D|}$, induced by the collection of maps $\hat{\iota}_{D_I}:D_I\to \tau_\epsilon^{\hat{X}}(D_I)$.

Now take $E\in\Spt(k)$. Applying $E$ to the diagram \eqref{eqn:CoMVDiag1} yields the diagram of presheaves on $|D|$
\[
I\mapsto i_{I*}[E(\tau_\epsilon^{\hat{X}}(D_I))].
\]
Applying $E$ to \eqref{eqn:CoMVDiag2} yields the diagram of presheaves on $|D|$
\[
I\mapsto i_{I*}[E(D_I)].
\]

\begin{Def} For $D\subset X$ a strict normal crossing subscheme and $E\in\Spt(k)$, we let $\tau_\epsilon^{\hat{X}}(D_\Zar)$ denote the diagram of co-presheaves \eqref{eqn:CoMVDiag1}, and set 
\[
E(\tau_\epsilon^{\hat{X}}(D)):=\holim_{I\neq\0}  i_{I*}[E(\iota_{I*}(\tau_\epsilon^{\hat{X}}(D_I))].
\]
We let $D_\Zar$ denote the diagram   \eqref{eqn:CoMVDiag2}, and set 
\[
E(D_\Zar):=\holim_{I\neq\0}  i_{I*}[E(D_{I\Zar})].
\]
\end{Def}
The natural transformation $\hat{\iota}_{|D|}$ yields the maps of presheaves on $|D|_\Zar$
\[
\hat{\iota}_{|D|}^*:E(\tau_\epsilon^{\hat{X}}(D))\to E(D_\Zar).
\]

\begin{prop}\label{prop:Retraction} The map $\hat{\iota}_{|D|}^*$ is a Zariski-local weak equivalence.
\end{prop}

\begin{proof} By Theorem~\ref{thm:main1}, the maps $\hat{\iota}_{D_I}$ are Zariski-local weak equivalences. The result follows from this and \cite{BousfieldKan}.
\end{proof} 

\begin{rem} One could also attempt a more direct definition of $\tau_\epsilon^{\hat{X}}(D)$ by just using our definition in the smooth case $i:W\to X$ and replacing the smooth $W$ with the normal crossing scheme $D$, in othere words, the co-presheaf on $|D|_\Zar$
\[
|D|\setminus F\mapsto \hat{\Delta}^*_{X\setminus F,|D|\setminus F}.
\]
Labeling this choice $\tau_\epsilon^{\hat{X}}(D)_{\na}$, and considering $\tau_\epsilon^{\hat{X}}(D)_{\na}$ as a constant diagram, we have the evident map of diagrams
\[
\phi:\tau_\epsilon^{\hat{X}}(D)\to \tau_\epsilon^{\hat{X}}(D)_{\na}
\]
We were unable to determine if $\phi$ induces a weak equivalence after evaulation on $E\in\Spt(k)$, even assuming that $E$ is homotopy invariant and satisfies Nisnevic excision. We were also unable to construct such an $E$ for which  $\phi$ fails to be a weak equivalence.  \end{rem}

\subsection{The punctured tubular neighborhood}

To define the punctured tubular neighborhood $\tau_\epsilon^{\hat{X}}(D)^0$, we proceed as follows: Fix an index $I\subset\{1,\ldots, m\}$, $I\neq\0$. Let $p:X'\to X$, $s:D_I\to X'$ be a Nisnevic neighborhood of $D_I$ in $X$, and let $|D|_{X'}=p^{-1}(|D|)$. Sending $X'\to X$ to $\Delta^n_{|D|_{X'}}$ gives us the pro-scheme $\hat{\Delta}^n_{X,|D|,D_I}$, and the closed embedding $\hat{\Delta}^n_{X,|D|,D_I}\to \hat{\Delta}^n_{X,D_I}$. Varying $n$, we have the cosimplicial pro-scheme  $\hat{\Delta}^*_{X,|D|,D_I}$, and the closed embedding $\hat{\Delta}^*_{X,|D|,D_I}\to \hat{\Delta}^*_{X,D_I}$.

For $U=D_I\setminus F$ an open subscheme of $D_I$, we set
\[
\tau_\epsilon^{\hat{X}}(D,D_I)^0(U):=\hat{\Delta}^*_{X\setminus F,U}\setminus \hat{\Delta}^*_{X\setminus F,|D|\setminus F,U}
\]
This forms the co-presheaf $\tau_\epsilon^{\hat{X}}(D,D_I)^0$ on $D_I$. The open immersions
\[
\hat{j}_I(U):\tau_\epsilon^{\hat{X}}(D,D_I)^0(U)\to \tau_\epsilon^{\hat{X}}(D,D_I)(U)
\]
give the map of co-presheaves
\[
\hat{j}_I:\tau_\epsilon^{\hat{X}}(D,D_I)^0\to \tau_\epsilon^{\hat{X}}(D,D_I).
\]

For $J\subset I$, we have the  map $\hat{i}_{J,I}: \hat{\Delta}^*_{X,D_I}\to \hat{\Delta}^*_{X,D_J}$ and
\[
\hat{i}_{J,I}^{-1}(\hat{\Delta}^*_{X,|D|,D_J})=\hat{\Delta}^*_{X,|D|,D_I}.
\]
Thus, we have the map $\hat{i}_{J,I}^0:\tau_\epsilon^{\hat{X}}(D,D_I)^0\to \tau_\epsilon^{\hat{X}}(D,D_J)^0$ and the diagram of co-presheaves on $|D|$
\begin{equation}\label{eqn:CoMVDiag3}
I\mapsto \iota_{I*}(\tau_\epsilon^{\hat{X}}(D_I)^0)
\end{equation}
which we denote by $\tau_\epsilon^{\hat{X}}(D)^0$.
The maps $\hat{j}_I$ define the map $\hat{j}_*:\tau_\epsilon^{\hat{X}}(D)^0\to
\tau_\epsilon^{\hat{X}}(D)$. Similarly, the maps $\pi_I:\tau_\epsilon^{\hat{X}}(D)\to X$ induce the map
\[
\pi^0:\tau_\epsilon^{\hat{X}}(D)^0\to X\setminus |D|,
\]
where we consider $X\setminus |D|$ the constant diagram.

\begin{Def} For $E\in\Spt(k)$, let $E(\tau_\epsilon^{\hat{X}}(D)^0)$ be the presheaf on $|D|_\Zar$,
\[
E(\tau_\epsilon^{\hat{X}}(D)^0):=\holim_{\0\neq I\subset \{1,\ldots, m\}}
\iota_{I*}E(\tau_\epsilon^{\hat{X}}(D_I)^0).
\]
\end{Def}

The map $\hat{j}_*$ defines the map of presheaves
\[
\hat{j}^*:E(\tau_\epsilon^{\hat{X}}(D))\to E(\tau_\epsilon^{\hat{X}}(D)^0).
\]
We let $E^{|D|}(\tau_\epsilon^{\hat{X}}(D))$ denote the homotopy fiber of $\hat{j}^*$. Via the commutative diagram
\[
\xymatrix{
E(X)\ar[d]_{\pi^*}\ar[r]^{j^*}&E(X\setminus |D|)\ar[d]^{\pi^{0*}}\\
E(\tau_\epsilon^{\hat{X}}(D))\ar[r]_{\hat{j}^*}&E(\tau_\epsilon^{\hat{X}}(D))^0
}
\]
we have the canonical map
\[
\pi^*_{|D|}:E^{|D|}(X)\to E^{|D|}(\tau_\epsilon^{\hat{X}}(D)).
\]

We want to show that the map $\pi^*_{|D|}$ is a weak equivalence, assuming that $E$ is homotopy invariant and satisfies Nisnevic excision.  We first consider a simpler situation.

We recall the simplicial structure on $E^{|D|}(\tau_\epsilon^{\hat{X}}(D))$ induced by the cosimplicial structure on the co-presheaves $\tau_\epsilon^{\hat{X}}(D_I)$ and $\tau_\epsilon^{\hat{X}}(D_I)^0$. In fact, letting
\begin{align*}
E(\tau_\epsilon^{\hat{X}}(D))_n&:=\holim_{I\neq\0}\iota_{I*}E(\hat{\Delta}^n_{\hat{X},D_{I\Zar}})\\
E(\tau_\epsilon^{\hat{X}}(D)^0)_n&:=\holim_{I\neq\0}\iota_{I*}E(\hat{\Delta}^n_{\hat{X},D_{I\Zar}}\setminus
\Delta^n_{|D|})\\
\end{align*}
and setting
\[
E^{|D|}(\tau_\epsilon^{\hat{X}}(D))_n:=\hofib(\hat{j}^*_n:
E(\tau_\epsilon^{\hat{X}}(D))_n\to E(\tau_\epsilon^{\hat{X}}(D)^0)_n),
\]
it follows from the fact that the homotopy limit is over a finite category that $E(\tau_\epsilon^{\hat{X}}(D))$,  $E(\tau_\epsilon^{\hat{X}}(D)^0)$ and $E^{|D|}(\tau_\epsilon^{\hat{X}}(D))$ are the total presheaves of spectra  associated to the simplical presheaves
\begin{align*}
&n\mapsto E(\tau_\epsilon^{\hat{X}}(D))_n\\
&n\mapsto E(\tau_\epsilon^{\hat{X}}(D)^0)_n\\
&n\mapsto E^{|D|}(\tau_\epsilon^{\hat{X}}(D))_n
\end{align*}
respectively. The map  $\pi^*_{|D|}$ is defined by considering $E^{|D|}(X)$ as a constant simplicial object. Let 
\[
\pi^*_{|D|,0}:E^{|D|}(X)\to E^{|D|}(\tau_\epsilon^{\hat{X}}(D))_0
\]
be the map of $E^{|D|}(X)$ to the 0-simplices of $E^{|D|}(\tau_\epsilon^{\hat{X}}(D))$.
  
\begin{lem} \label{lem:Diag}Suppose that $E$ satisfies Nisnevic excision and $D$ is a strict normal crossing subscheme of $X$. Then $\pi^*_{|D|,0}$ is a weak equivalence.
\end{lem}

Before we give the proof of this lemma, we prove an elementary result on Nisnevic neighborhoods.

\begin{lem} \label{lem:extension} Let $x$ be a smooth point on a finite type $k$-scheme $X$, let $\sO=\sO_{X,x}$ and let $Y=\Spec\sO$.  Let $D, E\subset Y$ be  smooth subschemes intersecting transversely and let $C=D\cap E$.
Let $p:\hat{Y}^h_D\to Y$, $q:\hat{Y}^h_C\to Y$ be the canonical maps and let $\hat{E}_D=p^{-1}(E)$, $\hat{E}_C=q^{-1}(E)$. Then there is a canonical $Y$-morphism $f:\hat{Y}^h_C\to \hat{Y}^h_D$. In addition,  $f^{-1}(\hat{E}_D)=\hat{E}_C$ and $f$ restricts to an isomorphism $\hat{E}_C\to \hat{E}_D$.
\end{lem}

\begin{proof} Since $Y$ is local, the pro-schemes $\hat{Y}^h_D$ and $\hat{Y}_C$ are represented by local $Y$-schemes. Since every Nisnevic neighborhood of $D$ in $Y$ gives a Nisnevic neighborhood of $C$ in $Y$, we have the canonical  $Y$-morphism $f:\hat{Y}^h_C\to \hat{Y}_D$. 

As $Y=\Spec\sO$ is local, we have a co-final family of   \'etale morphisms $Y'\to Y$ of the form $Y'=\Spec (\sO[T]/F)_G$, i.e., the localization of $\sO[T]/F$ with respect to some $G\in\sO[T]$, where $(\partial F/\partial T, F)$ is the unit ideal in $\sO[T]_G$. Those $Y'\to Y$ of this form which give a Nisnevic neighborhood of $D$ are those for which $F$ contains a linear factor, modulo the ideal $I_D$ of $D$. We have a similar description of a cofinal family of Nisnevic neighborhoods of $C$ in $Y$, and of $C$ in $E$. In particular, it follows that the Nisnevic neighborhoods of $C$ in $E$ which are the pull-back to $E$ of Nisnevic neighborhood of $D$ in $Y$ are co-final in the system of all Nisnevic neighborhoods of $C$ in $E$. Thus
\[
p^{-1}(E)\cong \hat{E}^h_C.
\]
A similar reasoning shows that $q^{-1}(C)\cong \hat{E}^h_C$, whence the lemma.
\end{proof}

\begin{proof}[proof of Lemma~\ref{lem:Diag}] Write $D$ as a sum, $D=\sum_{i=1}^mD_i$ with each $D_i$ smooth (but not necessarily irreducible), and with $m$ minimal. We proceed by induction on $m$.

For $m=1$, Nisnevic excision implies that the natural map
\[
E^{|D|}(X)\to E^{|D|}(\hat{X}^h_{|D|})
\]
is a weak equivalence. Noting that $ E^{|D|}(\hat{X}^h_{|D|})\to E^{|D|}(\tau_\epsilon^{\hat{X}}(D))_0$ is a weak equivalence since $D$ is smooth  proves the result in this case.

Before proceeding to the general case, we note that, for $W_1, W_2\subset Y$  closed subsets of $Y\in\Sm/k$, if $W=W_1\cup W_2$, the square
\[
\xymatrix{
E^{W_1\cap W_2}(Y)\ar[r]\ar[d]&E^{W_1}(Y)\ar[d]\\
E^{W_2}(Y)\ar[r]&E^W(Y)
}
\]
is homotopy cartesian. This follows easily from the Mayer-Vietoris property of $E$ with respect to Zariski open subschemes, which in turn is a direct consequence of Nisnevic excision.

Now assume the result for $m-1$. Let $D'=\sum_{i=2}^mD_i$ and write $|D|=D_1\cup |D'|$. We break up the diagram $\tau_\epsilon^{\hat{X}}(D)^0$ into three pieces: the single co-presheaf 
$\tau_\epsilon^{\hat{X}}(D_1)^0$, the diagram $\tau_\epsilon^{\hat{X}}(D_{*\ge2})^0$, which is $I\mapsto  \tau_\epsilon^{\hat{X}}(D_I)^0$, $1\not\in I$, and the diagram  $\tau_\epsilon^{\hat{X}}(D_{1,*\ge 2})^0$, which is $I\mapsto  \tau_\epsilon^{\hat{X}}(D_I)^0$, $1\in I$, $|I|\ge2$. We let $E(\tau_\epsilon^{\hat{X}}(D_{*\ge2})^0)$, $E(\tau_\epsilon^{\hat{X}}(D_{1,*\ge2})^0)$ denote the respective homotopy limits.

The evident maps of diagrams
\begin{align*}
&\tau_\epsilon^{\hat{X}}(D_{1,*\ge2})^0\to \tau_\epsilon^{\hat{X}}(D_1)^0\\
&\tau_\epsilon^{\hat{X}}(D_{1,*\ge2})^0\to \tau_\epsilon^{\hat{X}}(D_{*\ge2})^0
\end{align*}
give the homotopy cartersian diagram
\[
\xymatrix{
E(\tau_\epsilon^{\hat{X}}(D)^0)_0\ar[r]\ar[d]&E(\tau_\epsilon^{\hat{X}}(D_1)^0)_0\ar[d]\\
E(\tau_\epsilon^{\hat{X}}(D_{*\ge2})^0)_0\ar[r]&E(\tau_\epsilon^{\hat{X}}(D_{1,*\ge2})^0)_0
}
\]

If we perform the similar decomposition for the diagram $\tau_\epsilon^{\hat{X}}(D)$, giving the diagrams $\tau_\epsilon^{\hat{X}}(D_1)$, $\tau_\epsilon^{\hat{X}}(D_{*\ge2})$
and $\tau_\epsilon^{\hat{X}}(D_{1,*\ge2})^0$, we note that
\begin{align*}
&\tau_\epsilon^{\hat{X}}(D_{*\ge2})=\tau_\epsilon^{\hat{X}}(D')\\
&\tau_\epsilon^{\hat{X}}(D_{1,*\ge2})^0=\tau_\epsilon^{\hat{X}}(D'\cap D_1)
\end{align*}
so we have the homotopy cartesian square
\[
\xymatrix{
E(\tau_\epsilon^{\hat{X}}(D))_0\ar[r]\ar[d]&E(\tau_\epsilon^{\hat{X}}(D_1))_0\ar[d]\\
E(\tau_\epsilon^{\hat{X}}(D'))_0\ar[r]&E(\tau_\epsilon^{\hat{X}}(D'\cap D_1))_0
}
\]
Letting $E^{|D|_1}(\tau_\epsilon^{\hat{X}}(D_1))_0$, $E^{|D|'}(\tau_\epsilon^{\hat{X}}(D'))_0$ and $E^{|D|'_1}(\tau_\epsilon^{\hat{X}}(D'\cap D_1)_0$ be the homotopy fibers of the restriction maps $E(\tau_\epsilon^{\hat{X}}(D_1))_0\to E(\tau_\epsilon^{\hat{X}}(D_1)^0)_0$, etc., we have the homotopy cartesian square 
\[
\xymatrix{
E^{|D|}(\tau_\epsilon^{\hat{X}}(D))_0\ar[r]\ar[d]&E^{|D|_1}(\tau_\epsilon^{\hat{X}}(D_1))_0\ar[d]\\
E^{|D|'}(\tau_\epsilon^{\hat{X}}(D'))_0\ar[r]&E^{|D|'_1}(\tau_\epsilon^{\hat{X}}(D'\cap D_1))_0
}
\]

We now break up the terms in this last diagram using Mayer-Vietoris for closed subsets. Let $\pi_1:\hat{X}^h_{D_1}\to X$ be the canonical map, and let $|D'|_1=\pi_1^{-1}(|D'|)$. Then $|D|_1= |D_1|\cup |D'|_1$, and  $|D_1|\cap  |D'|_1 = |D_1|\cap |D'|$. 
By Nisnevic excision, the canonical maps
\begin{align*}
&E^{|D_1|\cap |D'|}(X)\to
E^{|D_1|\cap  |D'|}(\hat{X}^h_{D_1})=:E^{|D_1|\cap  |D'|}(\tau_\epsilon^{\hat{X}}(D_1))_0\\
&E^{|D_1|}(X)\to
E^{|D_1|}(\hat{X}^h_{D_1})=:E^{|D_1|}(\tau_\epsilon^{\hat{X}}(D_1))_0\\
\end{align*}
are weak equivalences, so   we have the homotopy cartesian square
\[
\xymatrix{
E^{|D_1|\cap |D'|_\Zar}(X)\ar[r]\ar[d]&
E^{|D_1|_\Zar}(X)\ar[d]\\
E^{|D'|_1}(\tau_\epsilon^{\hat{X}}(D_1))_0\ar[r]&
E^{|D|_1}(\tau_\epsilon^{\hat{X}}(D_1))_0
}
\]

Now let $|D_1|'_I$ be the inverse image of  $|D_1|$ under $\hat{X}^h_{D'_I}\to X$. Similarly, let 
$|D_1|_{1I}'$ and $|D'|_{1I}'$ be  the inverse images of $|D_1|$ and
$|D'|$ under $\hat{X}^h_{D'_I\cap D_1}\to X$, respectively. We define presheaves
\begin{align*}
&E^{|D_1|'}(\tau_\epsilon^{\hat{X}}(|D|'))_0(|D|\setminus F):=\holim_I
E^{|D_1|'_I\setminus F}(\hat{X}^h_{D'_I}\setminus F)\\
&E^{|D_1|_1'}(\tau_\epsilon^{\hat{X}}(D'\cap D_1))_0(|D|\setminus F):=
\holim_I E^{|D_1|_{1I}'\setminus F}(\hat{X}^h_{D'_I\cap D_1}\setminus F)\\
&E^{|D'|_1'}(\tau_\epsilon^{\hat{X}}(D'\cap D_1))_0(|D|\setminus F)
:=\holim_I E^{|D'|_{1I}'\setminus F}(\hat{X}^h_{D'_I\cap D_1}\setminus F)
\end{align*}
Using a similar argument as above, together with our induction hypothesis, we have homotopy cartesian squares
\[
\xymatrix{
E^{|D_1|\cap |D'|_\Zar}(X)\ar[r]\ar[d]&
E^{|D_1|'}(\tau_\epsilon^{\hat{X}}(|D|'))_0\ar[d]\\
E^{|D'|_\Zar}(X)\ar[r]&
E^{|D|'}(\tau_\epsilon^{\hat{X}}(|D|'))_0
}
\]
and
\[
\xymatrix{
E^{|D_1|\cap |D'|_\Zar}(X)\ar[r]\ar[d]&
E^{|D_1|_1'}(\tau_\epsilon^{\hat{X}}(D'\cap D_1))_0\ar[d]\\
E^{|D'|_1'}(\tau_\epsilon^{\hat{X}}(D'\cap D_1))_0\ar[r]&
E^{|D|_1'}(\tau_\epsilon^{\hat{X}}(D'\cap D_1))_0.
}
\]

We claim that the maps $\hat{X}^h_{D'_I\cap D_1}\to \hat{X}^h_{D'_I}$ and  $\hat{X}^h_{D'_I\cap D_1}\to \hat{X}^h_{D_1}$  induce Zariski-local  weak equivalences
\begin{align*}
&E^{|D'|_1}(\tau_\epsilon^{\hat{X}}(D_1))_0\to E^{|D'|_1'}(\tau_\epsilon^{\hat{X}}(D'\cap D_1))_0\\
&E^{|D_1|'}(\tau_\epsilon^{\hat{X}}(|D|'))_0\to E^{|D_1|'_1}(\tau_\epsilon^{\hat{X}}(D'\cap D_1))_0
\end{align*}

For the first map,  let $p:X'\to X$, $s:D_1\to X'$ be a Nisnevic neighborhood of $D_1$, and let $D''':=p^{-1}(D')$. By induction, the map
\[
E^{D''}(X')\to E^{D'}(\tau_\epsilon^{\hat{X'}}(D'')_0
\]
is a weak equivalence. Passing to the pro-scheme $\hat{X}^h_{D_1}$, we have the weak equivalence
\[
E^{\hat{D}'}(\hat{X}^h_{D_1})_0\to E^{\hat{D}'}(\tau_\epsilon^{\widehat{\hat{X}^h_{D_1}}}(\hat{D}'))_0,
\]
where $\hat{D}'$ is the inverse image of $D'$ under $\hat{X}^h_{D_1}\to
X$.  

The left-hand side of this equivalence is just $E^{|D'|_1}(\tau_\epsilon^{\hat{X}}(D_1))_0$. For each $I\subset \{2,\ldots, m\}$, we have the natural map of co-presheaves
\[
f:\hat{X}^h_{(D'\cap D_1)_I}\to (\widehat{\hat{X}^h_{D_1}})^h_{\hat{D}'_I}
\]
By Lemma~\ref{lem:extension}, for each $x\in (D'\cap D_1)_I$ the stalk $f_x$ of $f$ at $x$ defines a Nisnevic neighborhood of the localization $(\hat{D}'_I)_x$ of $\hat{D}'_I$ at $x$, and $f_x$ identifies the localization $(|D'|_1')_{I,x}$ of  $(|D'|_1')_I$ with  $(\hat{D}'_I)_x$. Thus, using Nisnevic excision for $E$, the map
\[
f^*:E^{\hat{D}'}(\tau_\epsilon^{\widehat{\hat{X}^h_{D_1}}}(\hat{D}'))_0\to
E^{|D'|_1'}(\tau_\epsilon^{\hat{X}}(D'\cap D_1))_0
\]
is a Zariski-local weak equivalence, as claimed.

The proof that the second map is a weak equivalence is similar, but simpler: one applies the above argument to each term in the diagram (indexed by the components of $D'$) defining $\tau_\epsilon^{\hat{X}}(|D|')_0$ to show that the map is a term-by-term weak equivalence; one need not use the induction hypothesis in this case.

Thus, putting the four homotopy cartesian squares together yields the homotopy cartesian square
\[
\xymatrix{
E^{|D_1|\cap |D'|_\Zar}(X)\ar[r]\ar[d]&
E^{|D_1|_\Zar}(X)\ar[d]\\
E^{|D'|_\Zar}(X)\ar[r]&
E^{|D|}(\tau_\epsilon^{\hat{X}}(D))
}
\]
Comparing this with the Mayer-Vietoris square
\[
\xymatrix{
E^{|D_1|\cap |D'|_\Zar}(X)\ar[r]\ar[d]&
E^{|D_1|_\Zar}(X)\ar[d]\\
E^{|D'|_\Zar}(X)\ar[r]&
E^{|D|_\Zar}(X)
}
\]
yields the result.
\end{proof}

\begin{prop}  \label{prop:Diag}Suppose that $E$ is homotopy invariant and satisfies Nisnevic excision, and $D$ is a strict normal crossing subscheme of $X$. Then $\pi^*_{|D|}$ is a weak equivalence.
\end{prop}

\begin{proof} Let $p:\Delta^n_{|D|}\to|D|$ be the projection. Applying Lemma~\ref{lem:Diag} to the strict normal crossing subscheme $\Delta^n_D$ of $\Delta^n_X$, the map
\[
p_*E^{|\Delta^*_D|_\Zar}(\Delta^*_X)\to E^{|D|}(\tau_\epsilon^{\hat{X}}(D))
\]
is a weak equivalence. Indeed,  $E^{|D|}(\tau_\epsilon^{\hat{X}}(D))$ is a simplicial object with $n$-simplices $p_*E^{|\Delta^n_D|}(\tau_\epsilon^{\widehat{\Delta^n_X}}(\Delta^n_D))_0$. Since $E$ is homotopy invariant, the map
\[
p^*:E^{|D|_\Zar}(X)\to p_*E^{|\Delta^*_D|_\Zar}(\Delta^*_X)
\]
is a (pointwise) weak equivalence, whence the result.
\end{proof}

We can now state and prove the main result for strict normal crossing schemes.

\begin{thm}\label{thm:main4} Let $D$ be a strict normal crossing scheme on some $X\in\Sm/k$ and take $E\in\Spt(k)$ which is homotopy invariant and satisfies Nisnevic excision. Then there is a natural distinguished triangle in $\sH\Spt(|D|_\Zar)$
\[
E^{|D|_\Zar}(X)\xrightarrow{\alpha_D} E(D_\Zar)\xrightarrow{\beta_D}E(\tau_\epsilon^{\hat{X}}(D)^0)
\]
\end{thm}
\begin{proof} By Proposition~\ref{prop:Diag}, we have a weak equivalence $E^{|D|_\Zar}(X)\to E^{|D|}(\tau_\epsilon^{\hat{X}}(D))$. By Proposition~\ref{prop:Retraction}, we have a Zariski local weak equivalence $E(\tau_\epsilon^{\hat{X}}(D))\to E(D_\Zar)$. Since $ E^{|D|}(\tau_\epsilon^{\hat{X}}(D))$ is by definition the homotopy fiber of the restriction map
 $E(\tau_\epsilon^{\hat{X}}(D))\to E(\tau_\epsilon^{\hat{X}}(D)^0)$, the result is proved.
 \end{proof}
 
 \begin{rem} From the work of Voevodsky \cite{VoevCross}, Ayoub \cite{Ayoub} and R\"ondigs \cite{Roendigs}, one has the general machinery of ``Grothendieck's six functors" in the setting of the motivic stable homotopy category $\SH(S)$ over a scheme $S$ (under certain conditions on $S$). For $i:D\to X$ a strict normal crossing scheme with complement $j:U\to X$, and for $\sE\in \SH(U)$, one has in particular the object $i^*j_*\sE$ of $\SH(S)$. One should view our construction $E(\tau_\epsilon^{\hat{X}}(D)^0)$ as a weak version of this, where first of all the input $E$ is in the category $\SH_{\A^1}(U)$, and the output $E(\tau_\epsilon^{\hat{X}}(D)^0)$ is in $\sH\Spt(|D|_\Zar)$. In particular,  $E(\tau_\epsilon^{\hat{X}}(D)^0)$ is only defined on Zariski open subsets of $|D|$, rather than on all of $\Sm/|D|$.
 \end{rem}

\section{Limit objects} The main interest in these constructions is to define the limit values $\lim_{t\to 0}E(X_t)$ for a semi-stable degeneration $\sX\to (C,0)$.

\subsection{Path spaces} Let $Y$ be a $k$-scheme. The free path space on $Y$ is defined to be the cosimplicial scheme $\sP_Y$ with $n$-cosimplices $Y^{n+2}$. The structure maps are defined as follows: Label the factors in $Y^{n+2}$ from 0 to $n+1$.  Send $\delta^n_i:[n]\to [n+1]$ to
\[
(y_0,\ldots y_{n+1})\mapsto (y_0,\ldots,y_{i-1}, y_i,y_i,y_{i+1}, y_{n+1})
\]
and send $s^n_i:[n]\to[n-1]$ to 
\[
(y_0,\ldots y_{n+1})\mapsto (y_0,\ldots,y_{i-1}, y_{i+1}, y_{n+1}).
\]
Projection on the first and last factors define the map $\pi:\sP_Y\to Y\times_kY$; we thus have two structures of a cosimplicial $Y$-scheme on $\sP_Y$: $\pi_1:\sP_Y\to Y$ and $\pi_2:\sP_Y\to Y$, with $\pi_i:=p_i\circ \pi$.

For a pointed $k$-scheme $(Y,y:\Spec k\to Y)$, set 
\[
\sP_Y(y):=\sP_Y\times_{(\pi_2,y)}\Spec k.
\]

Now let $p:\sY\to Y$ be a $Y$-scheme. Set
\[
\sP_{\sY/Y}(y):=\sY\times_{(p,\pi_1)}\sP_Y(y)
\]
We extend this definition to cosimplicial $Y$-schemes in the evident manner: if $\sY^\bullet\to Y$ is a cosimplicial $Y$-scheme, we have the bi-cosimplicial $Y$-scheme $\sP_{\sY^\bullet/Y}(y)$; the extension  to functors from some small category to cosimplicial $Y$-schemes is done in the same way.

Denoting the pointed $k$-scheme $(Y,y)$ by $Y_*$, we write $\sP_{Y_*}$ for $\sP_{Y}(y)$ and
$\sP_{\sY^\bullet/Y_*}$ for $\sP_{\sY^\bullet/Y}(y)$.

For $E\in\Spt(k)$, we have the simplicial spectrum $E(\sP_{\sY/Y_*}$.

\begin{lem}\label{lem:PathSpHom} Let  $\sY\to\A^1$ be a map of smooth $k$-schemes, and let $E\in\Spt(k)$ be homotopy invariant. Then the pull-back
\[
p_1^*:E(\sY)\to E(\sP_{\sY/\A^1}(0))
\]
is a weak equivalence.
\end{lem}

\begin{proof} $\sP_{\sY/\A^1}(0)$ has $n$-cosimplices $\sY\times_k\A^n$ and the map $p_1:\sP_{\sY/\A^1}(0)\to \sY$ on $n$-cosimplices is just the projection s $\sY\times_k\A^n\to\sY$. Since $E$ is homotopy invariant, 
\[
p_1^*:E(\sY)\to E(\sP_{\sY/\A^1}(0))
\]
is a termwise weak equivalence, hence a weak equivalence.
\end{proof}

\subsection{Limit structures} For our purposes, a semi-stable degeneration need not be proper, so even if this is somewhat non-standard terminology, we use the following definition:

\begin{Def} A {\em semi-stable degeneration} is a morphism $p:\sX\to (C,0)$, where
$(C,0)$ is a smooth pointed local curve over $k$, $C=\Spec\sO_{C,0}$, $\sX$ is a smooth irreducible $k$-scheme, $p$ is dominant,  smooth over $C\setminus 0$ and $p^{-1}(0):=X_0$ is a strict normal crossing divisor on $\sX$.
\end{Def}

For the rest of this section, we fix a  semi-stable degeneration $\sX\to(C,0)$. We denote the open complement of $X_0$ by $\sX^0$. We write $\G_m$ for the pointed $k$-scheme $(\A^1\setminus\{0\},1)$ (we sometimes tacitly forget the base-point).

Fix a uniformizing parameter $t\in\sO_{C,0}$, giving the morphism $t:(C,0)\to(\A^1_k,0)$, which restricts to $t:C\setminus0\to\G_m$. Let $p[t]:\sX\to\A^1$ be the composition $t\circ p$, and let $p[t]^0:\sX^0\to \G_m$ be the restriction of $p[t]$. Composing $p[t]$ with the canonical morphism $\tau_\epsilon^{\hat{\sX}}(X_0)^0\to
X^0$ yields the map 
\[
\hat{p}[t]^0:\tau_\epsilon^{\hat{\sX}}(X_0)^0\to \G_m.
\]
This makes  $\tau_\epsilon^{\hat{\sX}}(X_0)^0$ a diagram of co-presheaves (on $X_{0\Zar}$) of cosimplicial schemes over 
$\G_m$, so we have the diagram of cosimplicial co-presheaves on $X_{0\Zar}$:
\[
I\mapsto \sP_{\tau_\epsilon^{\hat{\sX}}(X_0)^0_I/\G_m}.
\]
Here $I$ runs over subsets of $\{1,\ldots, m\}$, where $X_0^1,\ldots, X_0^m$ are the irreducible components of $X_0$. We denote this diagram by
\[
\lim_{t\to0}X_t.
\]

Now let $E$ be in $\Spt(k)$. For each $I\subset\{1,\ldots, m\}$, we have the presheaf of bisimplicial spectra on $X_{0\Zar}$, $E(\sP_{\tau_\epsilon^{\hat{\sX}}(X_0)^0_I/\G_m})$, giving us the functor
\[
I\mapsto \tilde{E}(\sP_{\tau_\epsilon^{\hat{\sX}}(X_0)^0_I/\G_m}).
\]
where $\tilde{\ }$ means fibrant model. Taking the homotopy limit over $I$ of the associated diagram of presheaves of total spectra gives us the fibrant presheaf of spectra $E(\lim_{t\to0}X_t)$. Taking the global sections gives us the spectrum $E(\lim_{t\to0}X_t)(X_0)$, which we denote by $\lim_{t\to0}E(X_t)$.

\begin{prop} Suppose $E$ is homotopy invariant and satisfies Nisnevic excision. Then\\
(1) There is a canonical map  in $\sH\Spt(X_{0\Zar})$:
\[
E(X_{0\Zar})\xrightarrow{\gamma_{X}}E(\lim_{t\to0}X_t).
\]
(2)   If $X_0$ is smooth, then $\gamma_X$ is an isomorphism.
\end{prop}

\begin{proof} We repeat the construction of the diagram $\lim_{t\to0} X_t$, replacing the punctured tubular neighborhood $\tau_\epsilon^{\hat{\sX}}(X_0)^0$ with the full tubular neighborhood $\tau_\epsilon^{\hat{\sX}}(X_0)$. We let $E(\sP_{\tau_\epsilon^{\hat{\sX}}(X_0)/\A^1}(1))$ denote the homotopy limit of the diagram of presheaves on $X_{0\Zar}$
\[
I\mapsto \tilde{E}(\sP_{\tau_\epsilon^{\hat{\sX}}(X_0)_I/\A^1}(1)).
\]
By Lemma~\ref{lem:PathSpHom}, the map 
\[
p_1^*:E(\tau_\epsilon^{\hat{\sX}}(X_0))\to \tilde{E}(\sP_{\tau_\epsilon^{\hat{\sX}}(X_0)/\A^1}(1))
\]
is a weak equivalence. By Theorem~\ref{thm:main1} (and the homotopy invariance of $E$), 
\[
i_{X_0}^*:\tilde{E}(\tau_\epsilon^{\hat{\sX}}(X_0)) \to E(X_{0\Zar})
\]
is also a weak equivalence. We thereby have the canonical isomorphism
\[
 \tilde{E}(\sP_{\tau_\epsilon^{\hat{\sX}}(X_0)_I/\A^1}(1)) \cong E(X_0)
\]
The inclusions $\tau_\epsilon^{\hat{\sX}}(X_0)^0\to \tau_\epsilon^{\hat{\sX}}(X_0)$ give the map of diagrams
\[
\sP_{\tau_\epsilon^{\hat{\sX}}(X_0)^0_I/\A^1\setminus\{0\}}(1)\to
\sP_{\tau_\epsilon^{\hat{\sX}}(X_0)_I/\A^1}(1)
\]
and thus the restriction map
\[
\rho^*:E(\sP_{\tau_\epsilon^{\hat{\sX}}(X_0)/\A^1}(1))\to E(\lim_{t\to0}X_t).
\]
Composing  with the isomorphism described above yields the canonical map $\gamma_\sX$.

For (2),   fix a point $x\in X_0$. There is a Zariski neighborhood $U$ of $x$ in $X_0$ and a Nisnevic neighborhood $\sX'\to \sX$ of $U$ in $\sX$ which is isomorphic to a Nisnevic neighborhood of $U$ in $U\times\A^1$. Thus it suffices to prove the result in the case $\sX=X_0\times\A^1$, $(C,0)=(\A^1,0)$ and $p=p_2:\sX\to\A^1$.

For each smooth $k$-scheme $T$, the canonical map $\tau_\epsilon^{\hat{X_0\times\A^1}}(X_0\times 0)\to X_0\times\A^1$ induces a weak equivalence
\[
E(T\times_kX_{0\Zar}\times\A^1)\to E(T\times_k\tau_\epsilon^{\hat{X_0\times\A^1}}(X_0\times 0))
\]
Similarly, we have the weak equivalence
\[
E(T\times_kX_{0\Zar}\times\G_m)\to
E(T\times_k\tau_\epsilon^{\hat{X_0\times\A^1}}(X_0\times 0)^0).
\]
Applying these term-by-term with respect to the cosimplicial schemes  we have weak equivalences (assuming $\sX=X_0\times\A^1$)
\begin{align*}
&E(\sP_{\tau_\epsilon^{\hat{\sX}}(X_0)/\A^1}(1))(U)\to E(U\times\sP_{\A^1}(1))\\
&E(\sP_{\tau_\epsilon^{\hat{\sX}}(X_0)^0/\A^1\setminus\{0\}}(1))(U)\to E(U\times\sP_{\A^1\setminus\{0\}}(1))
\end{align*}
so we need only show that
\[
E(U\times\sP_{\A^1}(1))\to E(U\times\sP_{\A^1\setminus\{0\}}(1))
\]
is a weak equivalence for all $U$. 

As in the proof of (1), the projection $U\times\sP_{\A^1}(1)\to U$ induces a weak equivalence
\[
E(U)\to E(U\times\sP_{\A^1}(1))
\]
so we need only show that the projection $U\times\sP_{\G_m}\to U$ induces a weak equivalence
\[
E(U)\to E(U\times\sP_{\G_m})
\]

Let $\sk_NE(U\times\sP_{\A^1\setminus\{0\}}(1))$ denote the $N$-skeleton of the simplicial spectrum 
$E(U\times\sP_{\G_m})$, with respect to the simplicial structure
\[
n\mapsto E(U\times\sP_{\G_m})^n).
\]
Then, since $S^m$ is compact,
\begin{equation}\label{eqn:limit}
\pi_m(E(U\times\sP_{\G_m}))=\lim_{N\to\infty}
\pi_m(\sk_NE(U\times\sP_{\G_m})).
\end{equation}

We have the strongly convergent (homological) spectral sequence
\[
E^1_{p,q}=\pi_q\sk_NE(U\times\sP_{\G_m}^p)\Longrightarrow\pi_{p+q}\sk_NE(U\times\sP_{\G_m})
\]
Taking the limit in $N$ and using \eqref{eqn:limit}, this gives us the convergent spectral sequence
\[
E^1_{p,q}=\pi_qE(U\times\sP_{\G_m}^p)\Longrightarrow\pi_{p+q}E(U\times\sP_{\G_m})
\]

Let $T$ be a smooth $k$-scheme. The homotopy fiber sequence
\[
E^{T\times0}(T\times\A^1)\to E(T\times\A^1)\to E(T\times\G_m)
\]
is split by the map $i_1^*: E(T\times\G_m)\to E(T)$ and the weak equivalence
$i_1^*:E(T\times\A^1)\to E(T)$. Thus for every $n$, we have
\[
\pi_n(E(T\times\G_m))=\pi_n E(T)\oplus \pi_{n-1}E^{T\times0}(T\times\A^1).
\]

To make the calculation of the $E^2$-term easier, we replace the $E^1$-complex with the normalized subcomplex $NE^1$. Let $t_1,\ldots, t_n$ be the standard coordinates on $\G_m^n$, and define new coordinates $x_1,\ldots, x_n$ by
\[
x_i:=\begin{cases} t_it_{i+1}^{-1}&\text{ for } i=1,\ldots, n-1\\ t_n&\text{ for } i=n\end{cases}
\]
Under this change of coordinate, $\sP_{\G_m}$ becomes the cosimplicial scheme with $n$-cosimplices $\G_m^{n+1}$ and coface maps
\[
d_i(x_1,\ldots, x_n)= (x_1,\ldots, x_{i-1},1,x_i,\ldots, x_n)
\]
for $i=0,\ldots, n$. 
Using the notation $\Omega_{\P^1}E(T)$ for $E^{T\times0}(T\times\A^1)$, we thus have 
\[
NE^1_{p,q}=\pi_{q-p-1}(\Omega^{p+1}_{\P^1}E)(U)\oplus\pi_{q-p}(\Omega_{\P^1}^pE)(U).
\]
The differential $NE^1_{p,q}\to NE^1_{p-1,q}$ is the projection on $\pi_{q-p}(\Omega_{\P^1}^pE)(U)$ followed by the identity inclusion. Thus the spectral sequence degenerates at $E^2$ and gives
\[
\pi_nE(U\times\sP_{\A^1\setminus\{0\}}(1))=E^2_{0,n}=\pi_nE(U).
\]
The result follows easily from this.
\end{proof} 

 \section{The monodromy sequence} In this section, we construction the monodromy sequence for the limit object $E(\lim_{t\to0}X_t)$. As pointed out to us by J. Ayoub, one needs to restrict $E$ quite a bit. We give here a theory valid for presheaves of complexes of $\Q$-vector space on $\Sm/k$ which are homotopy invariant and satisfy Nisnevic excision.

 \subsection{Presheaves of complexes} For a noetherian ring $R$, we let $C_R$ denote the category of (unbounded) homological complexes of $R$-mod\-ules, $C_{R\ge0}$ the full subcategory of $C_R$ consisting of complexes which are zero in strictly negative degrees.
 
 By the Dold-Kan equivalence, we may identify $C_{R\ge0}$ with the category of simplicial $R$-modules $\Spc_R$. The forgetful functor $\Spc_R\to\Spc_*$ allows us to use the standard model structure on $\Spc_*$ to induce a model structure on $\Spc_R$, i.e., cofibrations are degreewise monomorphisms,  weak equivalences are homotopy equivalences on the geometric realization and fibrations are maps with the RLP for trivial cofibrations. This induces a model structure on $C_{R\ge0}$ with weak equivalence the quasi-isomorphisms; the suspension functor is the usual (homological) shift operator: $\Sigma C:=C[1]$, $C[1]_n:=C_{n-1}$, $d_{C[1],n}=-d_{C,n-1}$. This  model structure is extended to $C_R$ by identifying $C_R$ with the category of ``spectra in $C_{R\ge0}$", i.e., sequences $(C^0, C^1,\ldots)$ with bonding maps $\epsilon_n:C^n[1]\to C^{n+1}$. The model structure on $\Spt$ induces a model structure on spectra of simplicial $R$-modules, and thus a model structure on $C_R$, with weak equivalences the quasi-isomorphisms. Thus the homotopy category $\sH C_R$ is just the unbounded derived category $D_R$.
 
 Similarly, for a category $\sC$,  the model structure for the presheaf category $\Spt(\sC)$ gives a model structure for presheaves of complexes  on $\sC$, $C_R(\sC)$ with weak equivalences the pointwise quasi-isomorphisms, and homotopy category the derived category $D_R(\sC)$. We may introduce a topology (e.g., the Zariski or Nisnevic topology), giving the model categories $C_R(X_\Zar)$, $C_{R\Zar}(\Sm/B)$, $C_R(X_\Nis)$, $C_{R\Nis}(\Sm/B)$. These have homotopy categories equivalent to the  derived categories (on the small or big sites)  $D_R(X_\Zar)$,  $D_{R\Zar}(\Sm/S)$, $D_R(X_\Nis)$, $D_{R\Nis}(\Sm/S)$, respectively. Finally, we may consider the $\A^1$-localization, giving the Nisnevic-local $\A^1$-model structure $C_{R\A^1}(\Sm/B_\Nis)$ with homotopy category $D_{R\A^1}(B)$.
 
 Let $I$ be a category $F:I\to C_R$ a functor. Since we can consider $F$ as a spectrum valued-functor by various equivalences described above, we may form the complex $\holim_IF$. Explicitly, this is the following complex: One first forms the cosimplicial complex $\underline{\holim}_IF$ with $n$-cosimplices
 \[
 \underline{\holim}_IF^n:=\prod_{\sigma=(\sigma_0\to\ldots\to \sigma_n)\in \sN(I)_n} F(\sigma_n).
\]
For $g:[m]\to[n]$, with $g(m)=m'\le n$, the $\sigma$-component of the map $ \underline{\holim}_IF^n(g)$ sends $\prod x_\tau$ to $F(\sigma_{m'}\to \sigma_n)(x_{g^*(\sigma)})$. The complex $\holim_IF$ is then the total complex of the double complex $n\mapsto \underline{\holim}_IF^n$ and second differential the alternating sum of the coface maps. This construction being functorial and preserving quasi-iso\-morphisms, it passes to the derived category $D_R(\sC)$. If $I$ is a finite category, the construction commutes with filtered colimits, hence passes to the Zariski- and Nisnevic-local derived categories, as well as the $\A^1$-local versions.

 \begin{rems} (1) For a set $S$, let $R^S$ denote the free $R$-module on $S$. Sending a pointed space $(S,*)$ to the simplicial $R$-module $S_R$, with $S_{Rn}:=R^{S_n}/R^{\{*\}}$ defines the {\em $R$-localization functor} $\Spc_*\to \Spc_R$. This extends to the spectrum categories, and gives us the exact $R$-localization functor on homotopy category $\otimes R:\SH\to D_R$. The $R$-localization functor $\otimes R$ extends to all the model categories, in particular, we have the $R$-localization
 \[
 \otimes R:\SH_{\A^1}(B)\to D_{R\A^1}(B),
 \]

We can also take the $R$-localization of $\SH$ by performing a Bosufield localization, i.e., define $Z\in \Spt$ to be $\Q$-local if $\pi_n(Z)$ is a $\Q$-vector space for all $n$, and $E\to F$ a $\Q$ weak equivalence if $\Hom\_Spt(F,Z)\to \Hom_\Spt(E,Z)$ is an isomorphism for all $\Q$-local $Z$. Inverting the $\Q$-weak equivalences defines the $\Q$-local homotopy category $\SH_\Q$, and $\otimes\Q:\SH\to D_\Q$ identifies $\SH_\Q$ with $D_\Q$. This passes to the other homotopy categories we have defined, in particular, $\otimes \Q:\SH_{\A^1}(B)\to D_{\Q\A^1}(B)$ identifies $SH_{\A^1}(B)_\Q$ with 
$D_{\Q\A^1}(B)$.
\\
\\
(2) $D_{\Q\A^1}(k)$ is {\em not} the same as the ($\Q$-localized) big category of motives over $k$, $\DM(k)_\Q$.
\end{rems}

\subsection{The log complex}  Let $\sgn:S_n\to\{\pm1\}$ be the sign representation of the symmetric group $S_n$. Consider a presheaf  of $\Q$-vector spaces $E$ on $\Sm/k$. For $X, Y\in\Sm/k$, let $\alt_n:E(Y\times X^n)\to E(Y\times X^n)$ be the alternating projector
\[
\alt_n=\frac{1}{n!}\sum_{\sigma\in S_n}\sgn(\sigma)(\id_Y\times\sigma)^*,
\]
with $\sigma$ operating on $X^n$ by permuting the factors. 
let $E(Y,X^n)^\alt\subset E(Y\times X^n)$ be the image of $\alt_n$ and $E(Y,X^n)^\alt_\perp$ the kernel. 
We extend these constructions to presheaves of complexes $E$ by operating degreewise.

If $(X,*)$ is a pointed $k$-scheme, we have the inclusions
$i_j:  X^{n-1}\to Y\times X^n$ inserting the point $*$ in the $i$th factor.
For $E$ a presheaf of $\Q$-vector spaces, we let $E(Y\wedge X^{\wedge n})$  be the intersection of the kernels of the restriction maps
\[
(\id_Y\times i_j)^*:E(Y\times X^n)\to E(Y\times X^{n-1}).
\]
Letting $p_j:X^n\to X^{n-1}$ be the projection omitting the $j$th factor, the composition $(\id-p_n^*i_n^*)\circ\ldots\circ(\id-p_1^*i_1^*)$ gives a splitting 
\[
\pi_n:E(Y\times X^n)\to E(Y\wedge X^{\wedge n})
\]
 to the inclusion $E(Y\wedge X^{\wedge n})\to E(Y\times X^n)$.

Clearly $S_n$ acts on $E(Y\wedge X^{\wedge n})$ through its action on $X^n$; we let 
$E(Y, X^{\wedge n})^\alt$ and $E(Y, X^{\wedge n})^\alt_\perp$ be the image and kernel of $\alt_n$ on $E(Y\wedge X^{\wedge n})$. 
 
 Let $f:X\to \G_m$ be a morphism, $E$ a presheaf of $\Q$-vector spaces on $\Sm/k$.  Let $f_n:X\times\G_m^n\to X\times\G_m^{n+1}$ be the morphism
\[
f_n(x,t_1,\ldots, t_n):=(x,f(x), t_1,\ldots, t_n).
 \]
Denote the map $\alt_{n-1}\circ \pi_{n-1}\circ f_n^*:E(X,\G_m^{\wedge n})^\alt\to E(X,\G_m^{\wedge n-1})^\alt$ by
\[
\cup f:E(X,\G_m^{\wedge n})^\alt\to E(X,\G_m^{\wedge n-1})^\alt
\]
One checks that
\begin{lem} $(\cup f)^2=0$.
\end{lem}

\begin{proof} We work in the $\Q$-linear category $\Q\Sm/k$, with the same objects as $\Sm/k$, disjoint union being direct sum, and, for $X$, $Y$ connected, $\Hom_{\Q\Sm/k}(X,Y)$ is the $\Q$-vector space freely generated by the set $\Hom_{\Sm/k}(X,Y)$. Product over $k$ makes $\Q\Sm/k$ a tensor category. The map $\cup f$ is gotten by applying $E$ to the map $\cup f^\vee:X\times\G_m^{n-1}\to X\times\G_m^n$ in $\Q\Sm/k$:
\[
(x,t_1,\ldots, t_{n-1})\mapsto \alt[((x,f(x))-(x,1))\otimes t_1-1\otimes \ldots\otimes t_{n-1}-1)]
\]
and restricting to $E(X,\G_m^{\wedge n})^\alt$. But  $(\cup f^\vee)^2$ is
\begin{multline*}
(x,t_1,\ldots, t_{n-1})\mapsto\\ \ \alt[((x,f(x),f(x))-(x,1,f(x)) -(x,f(x),1)+(x,1,1))\otimes  \ldots\otimes t_{n-1}-1)
\end{multline*}
which is evidently the zero map.
\end{proof}

 Form the complex $E(\log_f)$ by 
 \[
 E(\log_f)_n:=E(X,\G_m^{\wedge n}))^\alt
 \]
with  differential $\cup f$. Since $E(\log_f)_0=E(X)$, we have the canonical map $\iota_X:E(X)\to E(\log_f)$.
 
We extend this definition to  an $I$-diagram of schemes over $\G_m$,  $f^\bullet :X^\bullet\to\G_m$ (with the $X^n\in\Sm/k$) by  
\[
E(\log_{f^\bullet}):=\holim_I(n\mapsto E(\log_{f^n}));
\]
similarly, we extend to $E$ a presheaf of complexes on $\Sm/k$ by taking the total complex of the double complex $n\mapsto E_n(\log_{f^\bullet})$. The map $\iota_X$ extends to
\[
\iota_{X^\bullet}:E(X^\bullet)\to E(\log_{f^\bullet}).
\]

 We consider as well a truncation of $E(\log_f)$. Recall that the stupid truncation $\sigma_{\ge n}C$ of a homological complex $C$ is the quotient complex of $C$ with 
 \[
 \sigma_{\ge n}C_m:=\begin{cases} C_m&\text{ for }m\ge n\\ 0&\text{ for }m<n.
 \end{cases}
 \]
 For $E$ a presheaf of abelian groups and $f:X\to \G_m$ a morphism in $\Sm/k$, set 
 \[
 E(\sigma_{\ge1}\log_f):=\sigma_{\ge1}E(\log_f).
 \]
We have the quotient map $N:E(\log_f)\to E(\sigma_{\ge1}\log_f)$, natural in $f$ and $E$.

We extend to $I$-diagrams $f^\bullet :X^\bullet\to\G_m$  and to presheaves of complexes as for $E(\log_f)$. The quotient map $N$ defined above extends to the natural map
\[
N:E(\log_{f^\bullet})\to E(\sigma_{\ge1}\log_{f^\bullet}).
\]
for $f^\bullet:X^\bullet\to \G_m$ an $I$-diagram of morphisms in $\Sm/k$, and $E\in C_\Q(k)$.

Finally, for $E\in C_\Q(k)$, define $E(-1)$ to be the presheaf of complexes
\[
E(-1)(X):=E(X\wedge\G_m)[1]:= \ker[E(X\times\G_m)\xrightarrow{i_1^*}E(X)][1].
\]

\begin{Def} \label{Def:alt} Let $E$ be be in $C_\Q(k)$. Call $E$ {\em alternating} if for every $X\in\Sm/k$ and every $n\ge0$, the alternating projection
\[
\alt_n:E(X\wedge\G_m^{\wedge n})\to E(X\wedge\G_m^{\wedge n})^\alt
\]
is a quasi-isomorphism.
\end{Def}

\begin{rems} (1) Clearly, $E$ is alternating if and only if $S_n$ acts via the sign representation on $H_mE(X\wedge\G_m^{\wedge n})$ for all $X$, $n$ and $m$. \\
\\
(2) Fix integers $1\le i\le n$. We have the split injection $\iota_{i,i+1}:
E(X\wedge\G_m^{\wedge n})\to E((X\times\G_m^{n-2})\wedge\G_m^{\wedge 2})$ by shuffling the $i, i+1$ coordinates to position $n-1, n$. In particular, we have the injection
\[
H_m(\iota_{i,i+1}):
H_mE(X\wedge\G_m^{\wedge n})\to H_mE((X\times\G_m^{n-2})\wedge\G_m^{\wedge 2}).
\]
Since $S_n$ is generated by simple transpositions, this shows that $E$ is alternating if and only if 
the exchange of factors in $\G_m\wedge\G_m$ acts by -1 on $H_mE(X\wedge\G_m\wedge\G_m)$ for all $X$ and $m$.\\
\\
(3) Suppose that $E\in C_\Q(k)$ is homotopy invariant and satisfies Nisnevic excision. Consider $\P^1$ as pointed by $\infty$. Then $E(X\wedge\P^1)$ is quasi-isomorphic to the suspension $E(X\wedge\G_m)[-1]$, hence $E$ is alternating if and only if the exchange of factors in $\P^1\wedge\P^1$ induces the identity on $H_mE(X\wedge\P^1\wedge\P^1)$ for all $X$ and $m$.

The homotopy invariance and Nisnevic excision properties of $E$ give a natural quasi-isomorphism of $E(X\wedge\P^1\wedge\P^1)$ with $E(X\wedge(\A^2/\A^2\setminus\{0\}))$, with the exchange of factors in $\P^1\wedge\P^1$ going over to the linear transformation $(x,y)\mapsto (y,x)$. If the characterstic of $k$ is different from 2, this transformation is conjugate to $(x,y)\mapsto (-x,y)$. Thus $E$ is alternating if and only if the map $[-1]:\P^1\to\P^1$, $[-1](x_0,x_1)=(x_0,-x_1)$, acts by the identity on $H_mE(X\wedge\P^1)$ for all $X$. \\
\\
(4) Call $E$ {\em oriented} if $E$ is an associative graded-commutive ring:
\[
\mu:E\otimes_\Q E\to E
\]
and (roughly speaking) $E$ admits a natural Chern class transformation
\[
c_1:\Pic\to H^2E
\]
satisfying the projective bundle formula: For $\sE\to X$ a rank $r$ vector bundle with associated projective space bundle $\P(\sE)\to X$ and tautological line bundle $\sO(1)$, $H^*E(\P(\sE))$ is a free $H^*E(X)$-module with basis $1, \xi, \ldots,\xi^{r-1}$, where $\xi=c_1(\sO(1))\in H^2E(\P(\sE))$. We do not assume that $c_1$ is a group homomorphism. The projective bundle formula and the fact that $[-1]^*\sO_{\P^1}(1)\cong\sO_{\P^1}(1)$ implies that an oriented $E$ is alternating. In particular, rational motivic cohomology, $\Q_\ell(*)$ \'etale cohomlogy, $\Q$-singular cohomology (with respect to a chosen embedding $k\to\C$) and rational algebraic cobordism $\MGL^{**}_\Q$ are all alternating.

On the other hand, rational motivic co-homotopy is alternating if $-1$ is a square in $k$, but is not alternating for $k=\R$. This is pointed out in \cite{Morel}: if $-1=i^2$, $[-1]$ is represented by the $2\times 2$ matrix with diagonal entries $i$ and $-i$. As this is a product of elementary matrices, one has an $\A^1$-homotopy connecting $[-1]$ and $\id$. To see the non-triviality of $[-1]$ for $k=\R$, let $[X,Y]$ denote the morphism $X\to Y$ in $\sH_{\A^1}\Spc_*(k)$. Morel defines a map (of sets) $[\P^1,\P^1]\to \phi(k)$, where $\phi(k)$ is the set of isomorphism classes of quadratic forms over $k$, and notes that the map $[u]$, 
\[
[u](x_0, x_1):=(x_0,ux_1), 
\]
goes to the class of the form $ux^2$. This map extends to a ring homomorphism
\[
\Hom_{\SH_{S^1}(k)}(\P^1,\P^1)\to \GW(k),
\]
where $\GW(k)$ is the Grothendieck-Witt ring (see also \cite[Lemma 3.2.4]{MorelLec2} for details).  Identifying $\GW(\R)$ with $\Z\times\Z$ by rank and signature, we see that $[-1]$ goes to the non-torsion element $(1,-1)$.

The example of motivic (co)homotopy is in fact universal for this phenomenon, so if $[-1]$ vanishes in $[\P^1,\P^1]$, then every $E\in C_\Q(k)$ satisfying homotopy invariance and Nisnevic excision is alternating.

We are grateful to F. Morel for explaining the computation of the transposition action on $\P^1\wedge\P^1$ in terms of quadratic forms and the Grothendieck-Witt group.

(4) Looking at the $\A^1$-stable homotopy category of $T$-spectra over $k$, $\SH(k)$, one can decompose the $\Q$-linearization $\SH(k)_\Q$ into the symmetric part $\SH(k)_+$ and alternating part $\SH(k)_-$ with respect to the exchange of factors on $\G_m\wedge\G_m$. Morel \cite{MorelMotCoh} has announced a result stating that $\SH(k)_-$ is in general equivalent to Voevodsky's big motivic category $\DM(k)_\Q$, and that $\SH(k)_+$ is zero if -1 is a sum of squares. This suggests that the alternating part $\SH_{S^1}(k)$ of the category of rational $S^1$-spectra $\SH_{S*1}(k)_\Q$ is closely related to the bigcategory of effective motives $\DM^\eff(k)$, but the exact relationship seems to be unclear at present.
\end{rems}

\begin{prop}\label{prop:Alt1} Let $E$ be in $C_\Q(k)$, $f:X\to \G_m$   an $I$-diagram of morphisms in $\Sm/k$.
\begin{enumerate}
\item The sequence 
\[
E(X)\xrightarrow{\iota_{X}} E(\log_f)\xrightarrow{N} E(\sigma_{\ge1}\log_f)
\]
identifies $E(\sigma_{\ge1}\log_f)$ with the quotient complex $E(\log_f)/E(X)$.
\item Suppose $E$ is alternating. Then there is a natural quasi-iso\-morphism $\alt:E(-1)(\log_f)\to E(\sigma_{\ge1}\log_f)$.
\end{enumerate}
\end{prop}

\begin{proof} It suffices to prove (1) for $E$ a presheaf of $\Q$-vector spaces, and $f:X\to \G_m$ a morphism in $\Sm/k$, where the assertion is obvious. Similarly, it suffices to construct a natural map $\theta_{E,X}:E(-1)(\log_f)\to E(\sigma_{\ge1}\log_f)$ for $E$ a presheaf of $\Q$-vector spaces, extend as above to a map in general,  and show that $\theta_{E,X}$ is a quasi-isomorphism for $E\in C_\Q(k)$ alternating and $f:X\to \G_m$ a morphism in $\Sm/k$.  

In fact,  for $E$ a presheaf of $\Q$-vector spaces and $n\ge1$, 
\[
E(-1)(\log_f)_n=\ker[(\id_X\times i)^*: E(X\times\G_m,\G_m^{\wedge n-1})^\alt\to 
E(X,\G_m^{\wedge n-1})^\alt]
\]
so $E(-1)(\log_f)_n$ is a subspace of $E(X,\G_m^{\wedge n})$; thus  $\alt_n$ defines a map $E(-1)(\log_f)_n\to E(\log_f)_n$. One easily checks that this defines a map of complexes
\[
\alt_*:E(-1)(\log_f)\to  E(\sigma_{\ge1}\log_f),
\]
as desired. 

Now suppose that $E$ is alternating, i.e., that 
\[
E(X,\G_m^{\wedge n})^\alt \to
E(X,\G_m^{\wedge n})
\]
 is a quasi-isomorphism 
for all $n$ and $X$.  This implies that the maps
\begin{align*}
&E(X\times\G_m,\G_m^{\wedge n-1})^\alt\to
E(X\times\G_m,\G_m^{\wedge n-1})\\
&E(X,\G_m^{\wedge n-1})^\alt\to
E(X,\G_m^{\wedge n-1})
\end{align*}
are quasi-isomorphisms, hence 
\[
\id_{X\wedge\G_m}\times\alt_{n-1}:
E(X\wedge\G_m,\G_m^{\wedge n-1})\to E(X\wedge\G_m,\G_m^{\wedge n-1})^\alt
\]
is a quasi-isomorphism. Since $E(X\wedge\G_m,\G_m^{\wedge n-1})=E(X,\G_m^{\wedge n})$, this implies that 
\[
\id_X\times\alt_n:E(X\wedge\G_m,\G_m^{\wedge n-1})^\alt\to E(X,\G_m^{\wedge n})^\alt
\]
is a quasi-isomorphism. As $\alt_*:E(-1)(\log_f)\to  E(\sigma_{\ge1}\log_f)$ is the map on the total complex of the double complexes 
\[
n\mapsto \id_X\times\alt_n:E(X\wedge\G_m,\G_m^{\wedge n-1})^\alt\to E(X,\G_m^{\wedge n})^\alt
\]
we see that $\alt_*$ is a quasi-isomorphism.
\end{proof}

\subsection{The log complex and path spaces} The monodromy sequence for $E(\lim_{t\to0}X_t)$ arises from the sequence of Proposition~\ref{prop:Alt1} by comparing the path space $E(\sP_{X/\G_m})$ and $E(\log_f)$ for $f:X\to \G_m$ a morphism in $\Sm/k$.

We use the Dold-Kan correspondence to rewrite $E(\sP_{X/\G_m})$ as a complex, namely: take for each $m$ the associated complex $E_m(sP_{X/\G_m})^*)$ of the simplicial abelian group $n\mapsto E_m(sP_{X/\G_m})^n)$. and then take the total complex of the double complex 
\[
m\mapsto E_m(sP_{X/\G_m})^*). 
\]
We write this complex as $E(\sP_{X/\G_m})$.

 We also have the normalized subcomplex $NE(\sP_{X/\G_m})$ of $E(\sP_{X/\G_m})$,  quasi-isomorphic to $E(\sP_{X/\G_m})$ via the inclusion. Recall that, for a simplicial abelian group $n\mapsto A_n$, the normalized complex $NA_*$ has
 \[
 NA_n:=\cap_{i=1}^n\ker d_i:A_n\to A_{n-1}
 \]
 with differential $d_0:NA_n\to NA_{n-1}$. We define $NE(\sP_{X/\G_m})$ by first taking the normalized subcomplex $NE_m(\sP_{X/\G_m})$ of $E_m(\sP_{X/\G_m}^*)$ for each $m$, and then forming the total complex of the double complex $m\mapsto NE_m(\sP_{X/\G_m})$.

  In particular, we have the inclusion of double complexes
 \[
 NE_*(\sP_{X/\G_m}^*)\subset E_*(\sP_{X/\G_m}^*);
 \]
 which gives for each $n$ the inclusion of single complexes
  \[
 NE_*(\sP_{X/\G_m}^n)\subset E_*(\sP_{X/\G_m}^n);
 \]
Recalling that $\sP_{X/\G_m}^n=X\times \G_m^n$, we thus have for each $n$ the inclusion of complexes
 \[
 NE_*(\sP_{X/\G_m}^n)\subset E_*(X\times \G_m^n),
 \]
 We may therefore apply the projections $\pi_n:E_*(X\times \G_m^n)\to
 E_*(X\times \G_m^{\wedge n})$ and $\alt_n$, giving the map
 \[
 \alt_n\circ\pi_n: NE_*(\sP_{X/\G_m}^n)\to E_*(X\times \G_m^{\wedge n})^\alt.
 \]
 
 \begin{lem} \label{lem:QIso} Suppose that $E$ is alternating. Then 
 \[
 \alt_n\circ\pi_n: NE_*(\sP_{X/\G_m}^n)\to E_*(X\times \G_m^{\wedge n})^\alt
 \]
  is a quasi-isomorphism.
 \end{lem}
 
 \begin{proof} The map $p_1^*:E(X)\to E(X\times\G_m)$ splits $i_1^*:E(X\times\G_m)\to E(X)$, so we have the natural splitting
 \[
 E(X\times\G_m)=E(X)\oplus E(X\wedge\G_m).
 \]
 Extending this to $E(X\times\G_m^n)$ by using the maps $i_j^*$ and $p_j^*$, we have the natural splitting
 \begin{equation}\label{eqn:split}
 E(X\times\G_m^n)=\oplus_{m=0}^n\oplus_{\substack{I\subset\{1,\ldots, n\}\\ |I|=m}}E(X\wedge\G_m^{\wedge I}).
 \end{equation}
 To explain the notation: For $I\subset \{1,\ldots, n\}$$E(X\wedge\G_m^{\wedge I})=E(X\wedge\G_m^{\wedge |I|})$, included in $E(X\times\G_m^n)$ by the composition
 \[
 E(X\wedge\G_m^{\wedge |I|})\subset E(X\times \G_m^{|I|})\xrightarrow{(\id_X\times p_I)^*}
 E(X\times\G_m^n)
 \]
 where $p_I:\G_m^n\to \G_m^{|I|}$ is the projection on the factors $i_1,\ldots, i_m$ if $I=\{i_1,\ldots, i_m\}$ with $i_1<\ldots i_m$.
 
 The action of $S_n$ on $E(X\times\G_m^n)$ preserves this decomposition, with $\sigma\in S_n$ mapping $E(X\wedge\G_m^{\wedge I})$ to $E(X\wedge\G_m^{\wedge \sigma^{-1}(I)})$ in the evident manner. 
 
Now, for a simplicial abelian group $A$, the inclusion $NA_n\to A_n$ is split by universal expressions in the face and degeneracy maps. If $n\mapsto C_{*n}$ is a simplicial complex, we can form the complex of normalized subgroups $N(C_{*_1,*_2})_{*_2=n}$ and take the homology $H_m(N(C_{*_1,*_2})_{*_2=n}, d_1)$, or we can  form the simplicial abelian group $n\mapsto H_m(C_{*n})$ and take the normalized subgroup $NH_m(C_{*_1, *_2},d_1)_{*_2=n}\subset H_m(C_{*_1,n},d_1)$. Using the universal spitting mentioned above, we see that the two are the same:
\[
H_m(N(C_{*_1,*_2})_{*_2=n}, d_1)=NH_m(C_{*_1, *_2},d_1)_{*_2=n}
\]
 Since $S_m$ acts by the sign representation on $H_mE(X\wedge\G_m^{\wedge m})$, it follows that, for $1\le j< m$, the diagonal embedding 
 \begin{align*}
& \delta_j:\G_m^{m-1}\to \G_m^m\\
&(t_1,\ldots, t_{m-1})\mapsto (t_1, \ldots t_j, t_j,t_{j+1},\ldots t_{m-1})
 \end{align*}
 induces the zero map on $H_mE(X\wedge\G_m^{\wedge m})$. Similarly, the inclusion 
  \begin{align*}
&i_n:\G_m^{n-1}\to \G_m^n\\
&(t_1,\ldots, t_{n-1})\mapsto (t_1,\ldots t_{m-1},1)
 \end{align*}
 is the zero map on $H_mE(X\wedge\G_m^I)$ if $n\in I$.
 
 From this, it is not hard to see that 
 \[
 NH_mE_*(NE_*(\sP_{X/\G_m}^n)=H_mE_*(X\wedge \G_m^{\wedge n}), 
 \]
 with respect to the decomposition of $E_*(\sP_{X/\G_m}^n)=E_*(X\times\G_m^n)$ given by \eqref{eqn:split}. Indeed, 
 \begin{align*}
 \ker(H_m(d_n))&=\ker(i_n^*:H_mE_*(X\times\G_m^n)\to H_mE_*(X\times\G_m^{n-1})\\
 &=\oplus_{\substack{I\subset\{1,\ldots, n\}\\n\in I}} H_mE_*(X\wedge\G_m^{\wedge I})
 \end{align*}
 It is then easy to show by descending induction on $i$ that 
 \[
 \cap_{j=i}^n\ker H_m(d_j)=\oplus_{\substack{I\subset\{1,\ldots, n\}\\\{i,\ldots, n\}\subset I}} H_mE_*(X\wedge\G_m^{\wedge I})
 \]
 from which our claim follows taking $i=1$. Thus the projection
 \[
 p_n:NE_*(\sP_{X/\G_m}^n)\to E_*(X\wedge\G_m^{\wedge n})
 \]
 is a quasi-isomorphism for each $n$. As $E$ is alternating, the alternating projection
 \[
 \alt_n:E_*(X\wedge\G_m^{\wedge n})\to E_*(X\wedge\G_m^{\wedge n})^\alt
 \]
 is a quasi-isomorphism as well, completing the proof. 
 \end{proof}
 
 \begin{lem}\label{lem:dif} Let $E$ be in $C_\Q(k)$. Let $\delta_0:E(X\times\G_m^{n})\to E(X\times\G_m^{n-1})$ be the map $[(\id_X,f),\id_{\G_m^{n-1}}]^*$. Then the diagram
 \[
 \xymatrix{
 E(X\times\G_m^{n})\ar[r]^{\delta_0}\ar[d]_{\alt_n\circ\pi_n}& E(X\times\G_m^{n-1})
 \ar[d]_{\alt_{n-1}\circ\pi_{n-1}}\\
 E(X\wedge\G_m^{\wedge n})^\alt\ar[r]_{\cup f}&E(X\wedge\G_m^{\wedge n-1})^\alt
 }
 \]
 commutes.
 \end{lem}
 
 \begin{proof} This follows directly from the definition of $\cup f$ and the fact that $\alt_n\circ\pi_n$ is the identity on $E(X\wedge\G_m^{\wedge n})^\alt$.
 \end{proof}
 
 \begin{prop}\label{prop:QIso} Let $E\in C_\Q(k)$ be alternating. Then the maps $ \alt_n\circ\pi_n: NE_*(\sP_{X/\G_m}^n)\to E_*(X\times \G_m^{\wedge n})^\alt$ define a quasi-isomorphism of total complexes
 \[
\alt\circ\pi: NE(\sP_{X/\G_m})\to E(\log_f)
\]
\end{prop}

\begin{proof} That the maps $\alt_n\circ\pi_n: NE_*(\sP_{X/\G_m}^n)\to E_*(X\times \G_m^{\wedge n})^\alt$ define a map of total complex $NE(\sP_{X/\G_m})\to E(\log_f)$ follows from Lemma~\ref{lem:dif} and the fact that, for each $n$, the differential $d_0$ on $NE(\sP_{X/\G_m}^n)$ is the restriction of  $\delta_0:E(X\times\G_m^{n})\to E(X\times\G_m^{n-1})$. Lemma~\ref{lem:QIso} implies that  $\alt\circ\pi$ is a quasi-isomorphism.
\end{proof}

We collect our results in 

\begin{thm}\label{thm:MonodromyDiag} Let $E\in C_\Q(k)$ be alternating, $f^\bullet :X^\bullet \to\G_m$ an $I$-diagram of morphisms in $\Sm/k$.  Consider the diagram
\[
\xymatrix{
&&E(\sP_{X/\G_m})&E(-1)(\sP_{X/\G_m})\\
&E(X^\bullet)\ar[r]^-{\iota_{X^\bullet}}\ar@{=}[dd]&NE(\sP_{X/\G_m})\ar[dd]_{\alt\circ\pi}\ar[u]^i&NE(-1)(\sP_{X/\G_m})\ar[d]_{\alt\circ\pi}\ar[u]^i\\
&&&E(-1)(\log_{f^\bullet})\\
0\ar[r]&E(X^\bullet)\ar[r]_-{\iota_{X^\bullet}}&E(\log_{f^\bullet})\ar[r]_-N&
E(\sigma_{\ge1}\log_{f^\bullet})\ar[u]_\alt\ar[r]&0
}
\]
Here the maps $i$ are the canonical inclusions and the maps $\iota_{X^\bullet}$ are the canonical maps 
given by the identities $E_n(\sP_{X^i/\G_m})_0=E_n(\log_{f^i})_0=E_n(X^i)$. Then
\begin{enumerate}
\item  The diagram commutes and   is natural in $E$ and $f^\bullet$.
\item All the maps in the diagram are maps of complexes.
\item All the vertical maps are quasi-isomorphisms
\item The bottom sequence is termwise exact.
\end{enumerate}
\end{thm}

\begin{proof} The first point follows by construction, the remaining assertions follow from the Dold-Kan correspondence, Proposition~\ref{prop:QIso} and Proposition~\ref{prop:Alt1}.
\end{proof}

\begin{cor} Let $W\in C_\Q(k)$ be alternating, $p:\sX\to (C,0)$ a semi-stable degeneration, $t\in \sO_{C,0}$ a uniformizing parameter. Then there is a
distinguished triangle in $D(X_{0\Zar})$
\[
E(\tau_\epsilon^{\hat\sX}(X_0)^0)\to E(\lim_{t\to0}X_t)\xrightarrow{N} E(-1)(\lim_{t\to0}X_t),
\]
natural in $(p,t)$ and in $E$.
\end{cor}

\begin{proof} The commutative diagram of Theorem~\ref{thm:MonodromyDiag} being natural in the choice of $I$-diagram and in $E$, one can extend the diagram directly to the case of a co-presheaf of cosimplicial $I$-diagrams 
\[
U\mapsto f^\bullet(U):X^\bullet(U)\to \G_m. 
\]
If we take $I$ to be finite, we can extend further to  a co-presheaf of cosimplicial $I$-diagrams $f^\bullet:X^\bullet\to \G_m$, with $X^\bullet(U)$ a  pro-scheme smooth over $k$,  and still preserve the quasi-isomorphisms and exactness. Feeding the $I$-diagram
\[
t\circ p:\tau_\epsilon^{\hat\sX}(X_0)^0\to\G_m
\]
to this machine and taking the distinguished triangle induced by the exact sequence of log complexes at the bottom of the diagram completes the proof.
\end{proof}

\begin{rem}\label{rem:MonodromySeq} If we splice together the long exact homotopy sequence for the monodromy distinguished triangle 
\[
E(\tau_\epsilon^{\hat\sX}(X_0)^0)\to E(\lim_{t\to0}X_t)\xrightarrow{N} E(-1)(\lim_{t\to0}X_t)
\]
with the localization distinguished triangle of Theorem~\ref{thm:main4} 
\[
E^{|D|_\Zar}(X)\xrightarrow{\alpha_D} E(D_\Zar)\xrightarrow{\beta_D}E(\tau_\epsilon^{\hat{X}}(D)^0)
\]
(both evaluated on $|D|=X_0$), we have the complex
\begin{multline}
\ldots\to E^{X_{0\Zar}}_n(X)\to E_n(X_{0\Zar})\to E_n(\lim_{t\to0}X_t)\xrightarrow{N}\\ E(-1)_n(\lim_{t\to0}X_t)\label{mult:MonoSeq}
\to  E^{X_{0\Zar}}_{n-2}(X)\to E_{n-2}(X_{0\Zar})\to \ldots
\end{multline}
If $k=\C$ and $E$ represents singular cohomology (for the classical topology) 
\[
E_n(Y)=H^{-n}(Y(\C),\Q)
\]
then Steenbrink's theorem \cite{Steenbrink} states that the above sequence is exact. The argument uses the mixed Hodge structure on all the terms together with a weight argument.

One should be able to define a natural geometric ``weight filtration" on $E(\lim_{t\to0}X_t)$ by using the stratification of $X_0$ by faces. However, for general $E$, this additional structure would probably not suffice to force the exactness of the above sequence. It would be interesting to give a general  additional structure on $E$ that would force the exactness.
\end{rem}

\section{Limit motives} We use our construction of limit cohomology, slightly modified, to give a definition of the limit motive of a semi-stable degeneration, as an object in the ``big" category of motives $\DM(k)$.

\subsection{The big category of motives} Voevodsky has defined the category of effective motives as the full subcategory $\DM^\eff_-(k)$ of the derived category of Nisnevic sheaves with transfer $D_-(\NST(k))$ consisting of those complexes with strictly homotopy invariant cohomology sheaves. 

In his thesis, Spitzweck \cite{SpitzweckThesis} defines a ``big" category of motives over a field $k$. R\"ondigs \cite{RoendigsBigMot} has also defined a big  category of motives over a noetherian base scheme $S$. To give the reader the main idea of both constructions, we quote from a recent letter from R\"ondigs \cite{RoendigsEmail}:

``One may construct a model category of simplicial pre\-sheaves with transfers on 
$\Sm/k$, in which the weak equivalences and fibrations are defined via the 
functor forgetting transfers. Via the Dold-Kan correspondence, there is an 
induced model structure on nonnegative chain complexes  of presheaves with 
transfers. Both may be stabilized with respect to $T$ or $\P^1$, in the sense of 
\cite{HoveySymSpec}.  The 
Dold-Kan correspondence extends accordingly. Since $T$ is a suspension already,
one can then pass to a model category of $\G_m$-spectra of integer-indexed chain 
complexes as well. For k a perfect field, results from \cite{FSV} show that the homotopy category of the 
latter model category contains Voevodsky's $\DM_{gm}$ as a full subcategory. "

We will use the $\P^1$-spectrum model.

\subsection{The cohomological motive}
We start with the category of pre\-sheaves with transfer $\PST(k)$ on $\Sm/k$, which is defined as in \cite{FSV} as the category of presheaves on the correspondence category $\Cor(k)$.
We let $C_{\ge0}(\PST(k))$ denote the model category of non-negative chain complexes in $\PST(k)$, with model structure induced from simplicial presheaves on $\Sm/k$, as described above. For $P\in C_{\ge0}(\PST(k))$, let $P(-1)$ denote the presheaf
\[
Y\mapsto \ker[P(Y\times\P^1)\xrightarrow{i_\infty^*} P(Y\times\infty)][2].
\]
where ``$\ker$' means the termwise kernel of the termwise split surjection $i_\infty^*$.  One has the adjoint  isomorphism
\[
\Hom_{C_{\ge0}(\PST(k))}(C\otimes\tilde\Z^{tr}_{\P^1}, C')\cong \Hom_{C_{\ge0}(\PST(k))}(C, C'(-1)[-2])
\]
so the bonding maps for $\P^1$-spectra in $C_{\ge0}(\PST(k))$ can  be just as well defined via maps
\[
C_n\to C_{n+1}(-1)[-2].
\]
We will use this normalization of the bonding morphisms from now on.

Now take $X\in\Sm/k$. For an integer $n\ge0$, we have the (homological) Friedlander-Suslin presheaf $\Z_{FS}^X(n)$:
\[
\Z_{FS}^X(n)(Y):=\Z_{FS}(n)(X\times Y):=C_{*-2n}(z_{q.fin}(\A^n))(X\times Y).
\]
We define 
\[
\delta_n:\Z_{FS}^X(n)\to \Z_{FS}^X(n+1)(-1)[-2]
\]
 by sending a cycle $W$ on $X\times Y\times\A^n$ to $W\times\Delta$, where $\Delta\subset \A^1\times\P^1$ is the graph of the inclusion $\A^1\subset \P^1$, and then reordering the factors to yield a cycle on $X\times Y\times\P^1\times\A^{n+1}$.
 
 \begin{Def} Let $X$ be in $\Sm/k$. The {\em cohomological motive} of $X$ is the sequence
 \[
\tilde{h}(X):=(\Z_{FS}^X(0), \Z_{FS}^X(1),\ldots, \Z_{FS}^X(n),\ldots)
 \]
 with the bonding morphisms $\delta_n$.
 \end{Def} 
 
 \begin{rem} One can also define the cohomological motive $h(X)\in \DM_{gm}(k)$ as the dual of the usual (homological) motive $m(X):=C^\Sus(\Z_\tr(X))$. For $X$ of dimension $d$, $h(X)(n)$ is actually in $\DM^\eff_-(k)$ for all $n\ge d$, and is represented by $\Z_{FS}^X(n)$. From this, one sees that the image of $\tilde{h}(X)$ in $\DM(k)$ is canonically isomorphic to $h(X)$. 
 
Also, one can  work in 
 $\DM^\eff_-(k)$ if one wants to define the cohomological motive of a diagram in $\Sm/k$ if the varieties involved have uniformly bounded dimension. Since our construction of limit cohomology uses varieties of arbitrarily large dimension, we need to work in $\DM(k)$.
 \end{rem}

 \subsection{The limit motive} It is now an easy matter to define the limit motive for a semi-stable degeneration. Let $\sX\to (C,0)$ be a semi-stable degeneration with parameter $t$ at 0; suppose the special fiber $X_0$ has irreducible components $X_0^1,\ldots, X_0^m$. We have the diagram of cosheaves on $X_{0\Zar}$, $\lim_{t\to0}X_t$, indexed by the non-empty subsets $I\subset \{1,\ldots, m\}$, which we write as
 \[
 I\mapsto [\lim_{t\to0}X_t]_I.
 \]
 Taking global sections on $X_0$ yields the diagram of cosimplicial schemes
  \[
 I\mapsto [\lim_{t\to0}X_t]_I(X_0).
 \]
 Applying $\tilde{h}$ gives us the digram of $\P^1$-spectra in $C_{\ge0}(\PST(k))$
 \[
 I\mapsto \tilde{h}([\lim_{t\to0}X_t]_I(X_0)).
 \]
 We then take the homotopy limit over this diagram forming the complex
 \[
\lim_{t\to0} \tilde{h}(X_t):=\holim_I \{ I\mapsto \tilde{h}( [\lim_{t\to0}X_t]_I(X_0))\}.
\]

\begin{Def} Let $\sX\to (C,0)$ be a semi-stable degeneration with parameter $t$ at 0. The limit cohomological motive $\lim_{t\to0}h(X_t)$ is the image of $\lim_{t\to0} \tilde{h}(X_t)$ in $\DM(k)$.
\end{Def}

Using the same procedure, we have, for $D\subset X$ a normal crossing scheme,  the motive of the tubular neighborhood $h(\tau_\epsilon^{\hat{X}}(D))$ and the  motive of the punctured tubular neighborhood $h(\tau_\epsilon^{\hat{X}}(D)^0)$. All the general results now apply for these cohomological motives. In particular, we have the monodromy distinguished triangle (for the $\Q$-motive)
\[
h(\tau_\epsilon^{\hat{\sX}}(X_0)^0)_\Q\to \lim_{t\to0}h(X_t)_\Q\to \lim_{t\to0}h(X_t)_\Q(-1)
\]
and the localization distinguished triangle
\[
h(X_0)\to h(\tau_\epsilon^{\hat{\sX}}(X_0)^0)\to h^{X_0}(X_0)
\]
From this latter triangle, we see that $h(\tau_\epsilon^{\hat{\sX}}(X_0)^0)$ is in $\DM_{gm}(k)$.

\section{Gluing  smooth curves}  We use the exponential map defined in \S\ref{sec:exp} to define an algebraic version of gluing  smooth curves along boundary components. We begin by recalling the construction of the moduli space of smooth curves with boundary components; for details we refer the reader to the article by Hain \cite{Hain}.

\subsection{Curves with boundary components}
For a $k$-scheme $Y$, a {\em smooth curve over $Y$} is a smooth proper morphism of finite type $p:\sC\to Y$ with geometrically irreducible fibers. We say that $\sC$ has genus $g$ if all the geometric fibers of are curves of genus $g$. A {\em boundary component} of $\sC\to Y$ consists of a section $x:Y\to\sC$ together with an isomorphism $v:\sO_Y\to x^*T_{\sC/Y}$, where $T_{\sC/Y}$ is the relative tangent bundle on $\sC$. Equivalently, $v$ is a nowhere vanishing section of $T_{\sC/Y}$ along $x$. A smooth curve with $n$ boundary components is $(\sC\to Y, (x_1,v_1),\ldots,(x_n,v_n))$ with all the $x_i$ disjoint. One has the evident notion of isomorphism of such tuples, so we can consider the functor $\sM_g^n$ on $\Sch_k$:
\begin{multline*}
\sM_g^n(Y):=\\ \{\text{smooth genus $g$ curves over $Y$ with $n$ boundary components} \}/\cong
\end{multline*}

For $n=0$, this is just the well-know functor of moduli of smooth curves, which admits the coarse moduli space $M_g$. For $n\ge1$, it is easy to show that a smooth curve over $Y$ with $n$ boundary components admits no automorphisms (over $Y$), from which it follows that $\sM_g^n$ is representable; we denote the representing scheme by $\sM_g^n$ as well.  

One can form a partial compactification of  $\sM_g^n$ by allowing {\em stable} curves with boundary components. As we will not require the full extent of this theory, we restrict ourselves to connected curves $C$ with a single singularity, this being an ordinary double point $p$.   We require that the boundary components are in the smooth locus of $C$. If $C$ is reducible, then $C$ has two irreducible components $C_1$, $C_2$; we also require that both $C_1$ and $C_2$ have at least one boundary component. As above, such data has no non-trivial automorphisms, which leads to the existence of a fine moduli space $\bar{\sM}^n_g$. We let $\sC_g^n\to \sM_g^n$ be the universal curve with universal boundary components $(x_1,v_1),\ldots, (x_n,v_n)$, and  $\bar{\sC}_g^n\to \bar{M}_g^n$  the extended universal curve.

The boundary $\del\bar{\sM}^n_g:=\bar{\sM}^n_g\setminus  {\sM}^n_g$ is a disjoint union of divisors
\[
\del\bar{\sM}^n_g:= D_{(g,n)}\coprod \coprod_{(g_1, g_2),(n_1,n_2)}D_{(g_1, g_2),(n_1,n_2)},
\]
where $D_{(g_1, g_2),(n_1,n_2)}$ consists of the curves $C_1\cup C_2$ with $g(C_i)=g_i$, and with $C_i$ having $n_i$ boundary components (we specify which component is $C_1$ by requiring $C_1$ to contain the boundary component $(x_1, v_1)$) and $D_{(g,n)}$ is the locus of irreducible singular curves. 

Let $(C, (x_1,v_1),\ldots)$ be a curve in $\del\bar{\sM}^n_g$ with singular point $p$. Let $C^N\to C$ be the normalization of $C$, and let $a,b\in C^N$ be the two points over $p$. By standard deformation theory, it follows that $\del\bar{\sM}^n_g$ is a smooth divisor in $\bar{\sM}^n_g$; let $N_{(g_1, g_2),(n_1,n_2)}$ denote the normal bundle of $D_{(g_1, g_2),(n_1,n_2)}$. Deformation theory gives a canonical identification  of the fiber of the punctured normal bundle  $N^0_{g_1, g_2, n_1, n_2}:=N_{(g_1, g_2),(n_1,n_2)}\setminus 0$ at $(C, (x_1,v_1),\ldots)$ with $\G_m$-torsor of isomorphisms
\[
T_{C^N,a}\otimes T_{C^N,b}\cong k(p).
\]

\subsection{Algebraic gluing}
We can now describe our algebraic construction of gluing curves. Fix integers $g_1, g_2$, $n_1, n_2\ge1$.  We defines the morphism
\[
\bar{\mu}:\sM_{g_1,n_1}\times\sM_{g_2,n_2}\to D_{g_1,g_2,n_1-1,n_2-1}.
\]
 by gluing $(C_1,(x_1,v_1),\ldots,(x_{n_1},v_{n_1}))$ and $(C_2,(y_1,w_1),\ldots, (y_{n_2},w_{n_2}))$ along $x_{n_1}$ and $y_1$, forming the curve $C:=C_1\cup C_2$ with boundary components $(x_1,v_1),\ldots,(x_{n_1-1},v_{n_1-1})$, $(y_2,w_2),\ldots, (y_{n_2},w_{n_2})$ and singular point $p$. We lift $\bar{\mu}$ to
\[
\mu:\sM_{g_1,n_1}\times\sM_{g_2,n_2}\to N^0_{g_1, g_2, n_1, n_2}
\]
using the isomorphism $T_{C_1,x_{n_1}}\otimes T_{C_2,y_1}\to k(p)$ which sends $v_{n-1}\otimes w_1$ to 1 and the identification of $(N^0_{g_1, g_2, n_1, n_2})_{C_1\cup C_2, \ldots}$ described above.

We now pass to the category $\SH_{\A^1}(k)$. Taking the infinite suspension, the map $\mu$ defines the map
\[
\Sigma^\infty\mu:\Sigma^\infty\sM_{g_1,n_1+}\wedge\Sigma^\infty\sM_{g_2,n_2+}\to \Sigma^\infty N^0_{g_1, g_2, n_1, n_2+}.
\]
Composing with our exponential map defined in \S\ref{sec:exp} gives us our gluing map
\[
\oplus:\Sigma^\infty\sM_{g_1,n_1+}\wedge\Sigma^\infty\sM_{g_2,n_2+}\to 
\Sigma^\infty\sM_{g_1+g_2,n_1+n_2-2+}.
\]

\begin{rems}(1)  If one fixes a curve $\sE:=(E,(x_1,v_1),(x_2,v_2))\in\sM_{1,2}$, one can form the tower under $\sE\oplus$
\[
\ldots\to\Sigma^\infty\sM_{g,n} \to \Sigma^\infty\sM_{g+1,n}\to\ldots, 
\]
and form the homotopy limit $\Sigma^\infty\sM_{\infty,n}$. If $E$ is an object of  $\SH_{\A^1}(k)$, one thus has the $E$-cohomology  $E^*(\sM_{\infty,n})$. For intance, this gives a possible definition of stable motivic cohomology or algebraic $K$-theory of smooth curves. However, it is not at all clear if this limit is independent of the choice of $\sE$. In the topological setting, one notes that the space $\sM_{1,2}(\C)$ is connected, so the limit cohomology, for example, is independent of the choice of $\sE$. On the contrary,  $\sM_{1,2}(\R)$ is not connected (the number of connected components in the real points of the curve corresponding to a real point of $\sM_{1,2}$ splits $\sM_{1,2}(\R)$ into disconnected pieces), so even there, the choice of $\sE$ plays a role.  It is also not clear if  $\sM_{\infty,n}$ is independent of $n$ (up to isomorphism in $\SH_{\A^1}(k)$).

\noindent
(2) In the topological setting,  the map $\oplus$ is the infinite suspension of a map
\[
\phi:\sM_{g_1,n_1}(\C)\times\sM_{g_2,n_2}(\C)\to 
\sM_{g_1+g_2,n_1+n_2-2}(\C),
\]
making $\amalg_{g,n}\sM_{g,n+2}(\C)$ into a topological monoid; the group completion is homotopy equivalent to the plus construction on the stable moduli space $\lim_{g\to\infty}\sM_{g,1}(\C)$ formed as in (1).   Letting $\sM_\infty(\C)^+$ denote this group completion, the group structure induces on $\Sigma^\infty\sM_\infty(\C)^+$ the structure of a Hopf algebra (this was pointed out to me by Fabian Morel), the co-algebra structure being the canonical one on a suspension spectrum. The functoriality of the exponential map $\exp^0$ as described in Remark~\ref{rem:expFunct} shows that the maps $\oplus$ make  $\bigvee_{g,n}\Sigma^\infty\sM_{g,n+2}$ into a biaglebra object in $\SH_{\A^1}(k)$. It is not clear if there is an analogous ``Hopf algebra completion" of $\bigvee_{g,n}\Sigma^\infty\sM_{g,n+2}$ in $\SH_{\A^1}(k)$.
\end{rems}

\section{Tangential base-points} Since motivic cohomology is represented in $\SH_{\A^1}(k)$, our methods are applicable to this theory. However, one can simplify the construction somewhat, since we are dealing with complexes of abelian groups rather than spectra. One can also achieve a refinement incorporating the multiplicative structure; this allows for a motivic definition of tangential base-points for the category of mixed Tate motives. 
 
\subsection{Cubical complexes} If we work with presheaves of complexes rather than presheaves of spectra, we can replace all our simplicial constructions with cubical versions. This enables any easy extension to the setting of differential graded algebras (\dga's), or even graded-commutative \dga's (\cdga's) if we work with complexes of $\Q$-vector spaces. We list the main results without proof here; the methods discussed in \cite[\S 2.5]{LevineHB} carry over without difficulty.

For a commutative ring $R$, we let $\bC_R$ denote the model category of (unbounded) complexes of $R$-modules and for a category $\sC$, let $\bC_R(\sC)$ denote the model category of presheaves of complexes on $\sC$. We copy the notation used for spectra, except that we denote as usual the homotopy category (derived category) by $\bD_R$, $\bD_R(\sC)$, respectively. For instance, for a scheme $X$, we have the model category of complexes of $R$-modules on the small Zariski site, $\bC(X_\Zar)$ and the derived category $\bD(X_\Zar)$. We have as well the model category of complexes of $R$-modules on the big Nisnevic site, $\bC(\Sm/S_\Nis)$, denoted $\bC_\Nis(S)$ and the derived category $\bD_\Nis(S)$.

The {\em cubical category} $\Cube$ has objects $\underline{n}$, $n=0,1,\ldots$. $\Cube$  is a subcategory of the category of finite sets, with $\underline{n}$ standing for the set $\{0,1\}^n$, with morphisms making $\Cube$ the smallest subcategory of finite sets containing the following  maps:
\begin{enumerate}
\item all inclusions $s_{i,n,\epsilon}:\{0,1\}^n\to \{0,1\}^{n+1}$, $\epsilon\in\{0,1\}$, $i=1,\ldots$, $n+1$, where $s_{i,n,\epsilon}$ is the inclusion inserting $\epsilon$ in the $i$th factor.
\item all projections $p_{i,n}:\{0,1\}^n\to \{0,1\}^{n-1}$, $i=1,\ldots, n$,  where  $p_{i,n}$ is the projection deleting the $i$th factor.
\item all maps $q_{i,n}:\{0,1\}^n\to \{0,1\}^{n-1}$, $i=1,\ldots, n-1$, $n\ge2$, defined by
\[
q_{i,n}(\epsilon_1,\ldots,\epsilon_n):=(\epsilon_1,\ldots,\epsilon_{i-1},\delta,\epsilon_{i+2},\ldots,\epsilon_n)
\]
with
\[
\delta:=\begin{cases}0&\text{ if }(\epsilon_i,\epsilon_{i+1}))=(0,0)\\
1&\text{ else.}
\end{cases}
\]
\end{enumerate}
A cubical object in a category $\sC$ is a functor $\Cube\to\sC$.

The basic cubical object in $\Sch$ is the sequence of $n$-cubes $\square^*:\Cube\to\Sm/k$. The operations of the projections $p_{i,n}$ and inclusions $s_{i,n}$ are the evident ones; $q_{i,n}$ acts by
\[
q_{i,n}(x_1,\ldots, x_n):=(x_1,\ldots, x_{i-1}, 1-(x_i-1)(x_{i+1}-1),x_{i+2},\ldots, x_n).
\]

Now let $P:\Cube\to \Mod_R$ be a cubical $R$-module. We have the {\em cubical realization} $|P|^c\in \bC_R$ with
\[
|P|^c_n:=P(\underline{n})/\sum_{i=1}^np_{i,n}^*(P(\underline{n-1})).
\]
The differential $d^c_n:|P|^c_n\to |P|^c_{n-1}$ is
\[
d^c_n:=\sum_{i=1}^n(-1)^is_{i,1}^*-\sum_{i=1}^n(-1)^is_{i,0}^*.
\]
$|-|^c$ is clearly a functor from the $R$-linear category of cubical $R$-modules to $\bC(R)$; in particular, if we apply $|-|^c$ to a complex of $R$-modules, we end up with a double complex. For a complex $C$,  also write $|C|^c$ for the total complex of this double complex, letting the context make the meaning clear.

\begin{ex}

For a presheaf of abelian groups $P$ on $\Sm/k$, we have cubical prescheaf $\sC^c(P)$ with
\[
\sC^c(P)(Y):=P(Y\times\square^*).
\]
Taking the cubical realization yields the  {\em cubical Suslin complex}  $C_*(P)^c$ with
\[
C_*(P)^c(Y):=|\sC^c(P)(Y)|.
\]
\end{ex}

The symmetric group $S_n$ acts on $C_n(P)^c$, we let $C_n(P)^c_{\alt}$ denote the subpresheaf of alternating sections.  One checks that the $C_n(P)^c_{\alt}$ form a subcomplex of $C_*(P)^c$.  If $P$ is a presheaf of $\Q$-vector spaces,  $C_n(P)^c_{\alt}$ is a canonical summand of  $C_n(P)^c$, with projection given by the idempotent $\Alt_n:=\frac{1}{n!}\sum_g\sgn(g)g$.

The main result on these constructions is
\begin{prop}\label{prop:Cubical} (1) There is a canonical homotopy equivalence of functors
\[
C_*\to C_*^c:\bC(k)\to \bC(k)
\]

\noindent
(2) If $P$ is a complex of presheaves of $\Q$-vector spaces, the inclusion
\[
C_*(P)^c_{\alt}\to C_*(P)^c
\]
is a quasi-isomorphism
\end{prop}

\begin{proof}[Sketch of proof; see \hbox{\cite[\S 5]{Loc}} for details.
] For (1), one uses the algebraic maps
\[
\square^n\to\Delta^n
\]
which collapse the faces $x_i=1$ to the vertex $(0,\ldots,0,1)$ to get a map $C_*\to C_*^c$. The homotopy inverse is given by triangulating the $\square^n$. 
For (2) one checks that $S_n$ acts by the sign representation on the homology sheaves of $C_*(P)^c$. The  projections $\Alt_n$ define a map of complexes $\Alt_*:C_*(P)^c\to C_*(P)^c_{\alt}$ which thus gives the inverse in homology.
\end{proof}

\subsection{Cubical tubular neightborhoods}  
For a closed embedding $i:W\to X$ in $\Sm/k$, set $\hat{\square}^n_{X, W}:=(\widehat{\square^n_X})^h_{\square^n_W}$, giving us the cubical pro-scheme
\[
\hat{\square}^*_{X,W}:\Cube\to\Pro\Sm/k
\]
We use the same notation for morphisms in the cubical setting as in the simplicial version, e.g., $\hat{i}_W:\square^*_W\to \hat{\square}^*_{X,W}$. We have as well the co-presheaf on $W_\Zar$ 
\[
\hat{\square}^n_{X, W_\Zar}(W\setminus F):=\hat{\square}^n_{X\setminus, W\setminus}
\] 
and the cubical co-presheaf 
\[
\tau_\epsilon^{\hat{X}}(W)^c:=\hat{\square}^*_{X, W_\Zar}.
\]
 
Now let $P$ be in $\bC(k)$. We define $P(\tau_\epsilon^{\hat{X}}(W)^c)_*$ to be the  complex of presheaves
\[
P(\tau_\epsilon^{\hat{X}}(W)^c)_*:=|P(\tau_\epsilon^{\hat{X}}(W)^c)|^c.
\]
We also have the alternating subcomplex $P(\tau_\epsilon^{\hat{X}}(W)^c)^\alt\subset P(\tau_\epsilon^{\hat{X}}(W)^c)$.

We have as well the punctured tubular neighborhood in cubical form
\[
\tau_\epsilon^{\hat{X}}(W)^{0c}:=\tau_\epsilon^{\hat{X}}(W)^c\setminus\square^*_{W_\Zar}
\]
on which we can evaluate $P$:
\[
P(\tau_\epsilon^{\hat{X}}(W)^{0c})_*:=|P(\tau_\epsilon^{\hat{X}}(W)^{0c})|^c.
\]
Let $P(\tau_\epsilon^{\hat{X}}(W)^{0c})^\alt\subset P(\tau_\epsilon^{\hat{X}}(W)^{0c})$ be the alternating subcomplex.

We let $EM:\bC\to \Spt$ be a choice of the Eilenberg-Maclane spectrum functor. Our main comparison result is
\begin{thm} (1) Let $i:W\to X$ be a closed embedding in $\Sm/k$.  For $P\in \bC(k)$, there are  natural isomorphisms in $\SH(W_\Zar)$
\begin{align*}
&EM(P(\tau_\epsilon^{\hat{X}}(W)^c))\cong EM(P)(\tau_\epsilon^{\hat{X}}(W))\\
&EM(P(\tau_\epsilon^{\hat{X}}(W)^{0c}))\cong EM(P)(\tau_\epsilon^{\hat{X}}(W)^0)
\end{align*}

\noindent
(2) If $P$ is a presheaf of complexes of $\Q$-vector spaces, then the inclusion 
\[
P(\tau_\epsilon^{\hat{X}}(W)^c)^\alt\to P(\tau_\epsilon^{\hat{X}}(W)^c)
\]
is a quasi-isomorphism.
\end{thm}

\begin{proof} Define
$P(\tau_\epsilon^{\hat{X}}(W))$ to be the total complex of the double complex associated to the simplicial complex $n\mapsto P(\tau_\epsilon^{\hat{X}}(W)^n)$. The   homotopy equivalence used in Proposition~\ref{prop:Cubical}(1) extends, by the functoriality of the Nisnevic neighborhood, to a homotopy equivalence
\[
P(\tau_\epsilon^{\hat{X}}(W))^c\sim P(\tau_\epsilon^{\hat{X}}(W))
\]
This yields a weak equivalence on the associated Eilenberg-Maclane spectra. Since the functor $EM$ passes to the homotopy category, we have a canonical isomorphism
\[
EM(P)(\tau_\epsilon^{\hat{X}}(W)))\cong EM(P(\tau_\epsilon^{\hat{X}}(W))).
\]
Putting these isomorphisms together completes the proof of the first assertion for the tubular neighborhood. The proof for the punctured tubular neighborhood is essentially the same. The second assertion follows from Proposition~\ref{prop:Cubical}(2).
\end{proof}

\subsection{The motivic \cdga} There are a number of different complexes which represent motivic cohomology; we will use the strictly functorial one of Friedlander-Suslin, $\Z_{FS}(q)$. Roughly speaking, one starts with the presheaf with transfers  of quasi-finite cycles $z_{q.fin}(\A^q)$, with value on $Y\in\Sm/k$ the cycles on $Y\times\A^q$ which are quasi-finite over $Y$, one forms the Suslin complex $C_*(z_{q.fin}(\A^q))$ and reindexes:
\[
\Z_{FS}(q)(Y)^n:=C_{2q-n}(z_{q.fin}(\A^q))(Y):=z_{q.fin}(\A^q)(Y\times\Delta^{2q-n}).
\]
 (see \cite[\S 2.4]{LevineHB} for a precise definition). This represents motivic cohomology Zariski-locally:
 \[
 H^p(X,\Z(q))=\H^p(X_\Zar,\Z_{FS}(q)).
 \]
 
We will use the cubical version $\Z_{FS}(q)^c$.
\[
\Z_{FS}(q)^{c, n}:=C_{2q-n}(z_{q.fin}(\A^q))^c(Y).
\]
By Proposition~\ref{prop:Cubical}, $\Z_{FS}(q)^c$ is quasi-isomorphic to $\Z_{FS}(q)$. 

Passing to $\Q$-coefficients, we have the quasi-isomorphic alternating subcomplex $\Q_{FS}(q)^c_\alt\subset
\Q_{FS}(q)^c$. We may also symmetrize with respect to the coordinates in the $\A^q$ in $z_{q.fin}(\A^q)$; it is shown in \cite{LevineHB} that the inclusion
\[
\Q_{FS}(q)^c_{\alt,\sym}\subset
\Q_{FS}(q)^c_\alt
\]
is also a quasi-isomorphism.

The product map
\[
z_{q.fin}(\A^q)(\square^n\times Y)\otimes z_{q.fin}(\A^{q'})(\square^{n'}\times Y)\to
z_{q.fin}(\A^{q+q'})(\square^{n+n'}\times Y)
\]
makes the graded complex 
\[
\tilde{\sN}_\Z:=\oplus_{q\ge0}\Z_{FS}(q)^c
\]
into a presheaf of Adams-graded \dga's on $\Sm/k$ (with Adams grading $q$). Passing to $\Q$-coefficients, and following the product with the alternating and symmetric projections makes 
\[
\sN:=\oplus_{q\ge0}\Q_{FS}(q)^c_{\alt,\sym}
\]
a presheaf of Adams-graded \cdga's,  the {\em motivic } \cdga\ on $\Sm/k$.

We let $\sN\to\sN^\fib$ denote a fibrant model of $\sN$ in the model category of (Adams-graded) \cdga's on $\Sm/k$, where the weak equivalences are Adams-graded quasi-isomorphisms of \cdga's for the Zariski topology. 

\begin{rems} (1)  Since $\sN$ is strictly homotopy invariant \cite[Theorem 4.2]{FSV}, $\sN^\fib$ is homotopy invariant.
\\
(2) In case $k$ admits resolution of singularities (i.e., $\Char k=0$) the canonical map $\Z_{FS}(q)\to \Z_{FS}(q)^\fib$ is a pointwise weak equivalence \cite[Theorem 7.4]{FSV}. Thus, in this case, we can use $\sN$ instead of $\sN^\fib$.
\end{rems}

\subsection{The specialization map}
We consider the situation of a smooth curve $C$ over our base-field $k$ with a $k$-point x. We let $\sO$ denote the local ring of $x$ in $C$, $K$ the quotient field of $\sO$ and choose a uniformizing parameter $t$, which we view as giving a map
\[
t:\Spec\sO\to \A^1.
\]
sending $x$ to 0.

Letting $i_x:x\to \Spec \sO$ be the inclusion, we have the restriction map
\[
i_x^*:\sN(\sO)\to \sN(k(x)),
\]
which is a morphism of Adams-graded \cdga's. In this section, we extend $i_x^*$ to a map
\[
sp_t:\sN(K)\to \sN(k(x))
\]
in the derived category of Adams-graded \cdga's over $\Q$ (denoted $\sH(\text{\cdga}_\Q)$).

First, if we apply $\sN$ to $\square^*\times Y$ and take the alternating projection again, we have the presheaf of \cdga's $\sN(\square^*_\alt)$ and the quasi-isomorphism of presheaves of \cdga's
\[
\iota:\sN\to \sN(\square^*_\alt).
\]

Next, write $\hat{\square}^{m0}_{C,x}$ for $\hat{\square}^m_{C,x}\setminus \square^m_x$, and consider the cubical punctured tubular neighborhood $\Z_{FS}(q)^c(\tau_\epsilon^{\hat{C}}(x)^{0c})$.  The product map
\begin{multline*}
z_{q.fin}(\A^q)(\square^n\times \hat{\square}^{m0}_{C,x})\otimes 
z_{q.fin}(\A^{q'})(\square^{n'}\times \times \hat{\square}^{m'0}_{C,x}))\\\to
z_{q.fin}(\A^{q+q'})(\square^{n+n'}\times \times \hat{\square}^{m+m'0}_{C,x}))
\end{multline*}
makes  $\oplus_{q\ge0}\Z_{FS}(q)^c(\tau_\epsilon^{\hat{\Spec\sO}}(x)^{0c})$ into an Adams-graded \dga; taking the alternating projection in both the $\square^n$ and $\hat{\square}^{m0}_{C,x}$ variables, and the symmetric projection in $\A^q$ and applying the fibrant model gives a presheaf of Adams-graded \cdga's, denoted $\sN^\fib(\tau_\epsilon^{\hat{\Spec\sO}}(x)^{0c}_\alt)$.  

Similarly, we perform this construction using the full tubular neighborhood, giving the presheaf  $\sN^\fib(\tau_\epsilon^{\hat{C}}(x)^c_\alt)$, and the commutative diagram of Adams-graded \cdga's:
\[
\xymatrix{
\sN(k(x))\ar[d]_\iota&\sN(\sO)\ar[r]^\res\ar[d]_{\pi_\sO^*}\ar[l]_{i_x^*}&\sN(K)\ar[d]^\iota\ar[d]^{\pi^*_K}\\
\sN^\fib(\square^*_\alt)(k(x))&\sN^\fib(\tau_\epsilon^{\hat{C}}(x)^{c}_\alt)\ar[r]_-\res\ar[l]^{i_x^*}
&\sN^\fib(\tau_\epsilon^{\hat{C}}(x)^{0c}_\alt)
}
\]

Replacing $(C,x)$ with $(\A^1,0)$ and using $\A^1$ and $\G_m$ instead of $\Spec \sO$ and $\Spec K$ yields the commutative diagram of Adams-graded \cdga's
\[
\xymatrix{
\sN(k(0))\ar[d]_\iota&\sN^\fib(\A^1)\ar[r]^-\res\ar[d]_{\pi^*_{\A^1}}\ar[l]_-{i_0^*}
&\sN^\fib(\G_m)\ar[d]^\iota\ar[d]^{\pi^*_{\G_m}}\\
\sN^\fib(\square^*_\alt)(k(x))&\sN^\fib(\tau_\epsilon^{\hat{\A^1}}(0)^{c}_\alt)\ar[r]_-\res\ar[l]^-{i_0^*}
&\sN^\fib(\tau_\epsilon^{\hat{\A^1}}(0)^{0c}_\alt).
}
\]
By Corollary~\ref{cor:Retract} and Corollary~\ref{cor:DistTriang2}, the maps $\pi^*_{\A^1}$ and $\pi^*_{\G_m}$  are quasi-isomorphisms of complexes, hence  quasi-isomorphisms of Adams-graded \cdga's. Since $\sN^\fib$ is homotopy invariant, the maps $\iota$ are quasi-isomorphisms of Adams-graded \cdga's. 

Finally, the map $t$ induces the commutative diagram of Adams-graded \cdga's
\[ 
\xymatrix{
\sN(k(x))&\sN^\fib(\tau_\epsilon^{\hat{C}}(x)^{c}_\alt)\ar[r]^-\res\ar[l]_-{i_x^*}
&\sN^\fib(\tau_\epsilon^{\hat{C}}(x)^{0c}_\alt) \\
\sN(k(0))\ar[u]^{t^*}&\sN^\fib(\tau_\epsilon^{\hat{\A^1}}(x)^{c}_\alt)\ar[r]_-\res\ar[u]^{t^*}\ar[l]^-{i_0^*}
&\sN^\fib(\tau_\epsilon^{\hat{\A^1}}(x)^{0c}_\alt)\ar[u]_{t^*}.
}
\]
Since $t:(C,x)\to (\A^1,0)$ is a Nisnevic neighborhood of $0$ in $\A^1$, all three maps $t^*$ are  isomorphisms. Putting these diagrams together and inverting the quasi-isomorphisms $\iota$, $t^*$, $\pi^*_{\A^1}$ and $\pi^*_{\G_m}$ yields the commutative diagram in $\sH(\cdga_\Q)$:
\begin{equation}\label{eqn:ComDiag1}
\xymatrix{
\sN(k(x))&\sN(\sO)\ar[l]_{i_x^*}\ar[r]^-\res\ar[d]_{\phi_\sO^*}&\sN(K)\ar[d]^{\phi^*_K}\\
\sN(k(0))\ar[u]^{t^*}&\sN^\fib(\A^1)\ar[r]_-\res\ar[l]^-{i_0^*}
&\sN^\fib(\G_m)
}
\end{equation}

\begin{Def} Let $i_1:\Spec k\to \G_m$ be the inclusion. The map $sp_t:\sN(K)\to \sN(k(x))$ in $\sH(\cdga_\Q)$ is defined to be the composition
\[
\sN(K)\xrightarrow{\phi^*_K}\sN^\fib(\G_m)\xrightarrow{i_1^*} \sN^\fib(k)\cong \sN(k)=\sN(k(0))\xrightarrow{t^*}
\sN(k(x)).
\]
\end{Def}

\begin{prop}\label{prop:Com1} The diagram in  $\sH(\cdga_\Q)$ 
\[
\xymatrix{
\sN(\sO)\ar[r]^-\res\ar[rd]_{i_x^*}&\sN(K)\ar[d]^{sp_t}\\
&\sN(k(x))
}
\]
commutes.
\end{prop}

\begin{proof} Since $\sN^\fib$ is homotopy invariant, the maps
\[
i_0^*, i_1^*:\sN^\fib(\A^1)\to \sN(k)
\]
are equal in  $\sH(\cdga_\Q)$. The proposition follows directly from this and a  chase of the commutative diagrams defined above.
\end{proof}

\subsection{The specialization functor} For a field $k$, we have the triangulated category $\DTM(k)$ of {\em mixed Tate motives} over $k$, this being the full triangulated subcategory of Voevodsky's triangulated category of motives (with $\Q$-coefficients), $\DM_{gm}(k)_\Q$, generated by the Tate objects $\Q(n)$, $n\in\Z$. 

We will also use in this section the derived category of {\em cell modules} over an Adams-graded \cdga\ $\sA$,  $\DCM(\sA)$. This construction was introduced in \cite{KrizMay}; we refer the reader to the discussion in \cite[\S 5]{LevineHB} for the properties of $\DCM$ we will be using below.

Let $\sO$ be as in the previous section the local ring of a $k$-point $x$ on a smooth curve $C$ over $k$, with quotient field $K$. The map $sp_t:\sN(K)\to \sN(k(x))$ yields an exact tensor functor
\[
sp_t:\DTM(K)\to \DTM(k(x))
\]

Indeed, as discussed in \cite[\S5.5]{LevineHB}, Spitzweck's representation theorem gives a natural equivalence of $\DTM(k)$ with the derived category $\DCM(\sN(k))$ of cell modules over the Adams-graded \cdga\ $\sN(k)$, as triangulated tensor $\Q$-tensor categories.

The functor  $\DCM$ associating to an Adams-graded $\Q$-\cdga\ $\sA$  the triangulated tensor category $\DCM(\sA)$ takes quasi-isomorphisms to  triangulated tensor  equivalences, hence  $\DCM$ descends to a well-defined pseudo-functor on $\sH(\cdga_\Q)$.
Thus, we may make the following
\begin{Def} Let $\sO$ be  the local ring of a $k$-point $x$ on a smooth curve $C$ over $k$, with quotient field $K$ and uniformizing parameter $t$.  Let $sp_t:\DTM(K)\to \DTM(k(x))$ be the exact tensor functor induced by $\DCM(sp_t):\DCM(\sN(K))\to\DCM(\sN(k(x))$, using  Spitzweck's representation theorem to identify the derived categories of cell modules with the appropriate category of mixed Tate motives.
\end{Def}

\begin{rem} (1) The discussion in \cite[\S5.5]{LevineHB}, in particular, the statement and proof of Spitzweck's representation theorem, is in the setting of motives over a field. However, we now have available a reasonable triangulated category $\DM(S)$ of motives over an arbitrary base-scheme $S$ (see \cite{VoevSimplMotives}), and we can thus define the triangulated category of mixed Tate motives over $S$, $\DTM(S)$, as in the case of a field. 

Furthermore, if $S$ is in $\Sm/k$ for $k$ a field of characteristic zero, then $\sN(S)$  has the correct cohomology, i.e.
\[
H^n(\sN(S))=\oplus_{q\ge0}H^n(S,\Q(q)),
\]
and one has the isomorphism
\[
H^n(S,\Z(q))\cong \Hom_{\DM(S)}(\Z,\Z(q)).
\]
This is all that is required for the argument in  \cite[\S5.5]{LevineHB} to go through, yielding the equivalence of the triangulated tensor category of cell modules $\DCM(\sN(S))$ with  $\DTM(S)$.\\
\\
(2) Joshua \cite{Joshua} has {\em defined} the triangulated category of $\Q$ mixed Tate motives over $S$ as $\DCM(\sN(S))$; the discussion in (1) shows that this agrees with the definition as a subcategory of $\DM(S)_\Q$.
\end{rem}

With these remarks, we can now state the main compatibility property of the functor $sp_t:\DTM(K)\to \DTM(k(x))$.

\begin{prop} \label{prop:Com2} Let $\sO$ be  the local ring of a $k$-point $x$ on a smooth curve $C$ over $k$, with quotient field $K$ and uniformizing parameter $t$. Let $i_x^*:\DTM(\sO)\to\DTM(k)$ and $j^*:\DTM(\sO)\to\DTM(K)$ be the functors induced by the inclusions $i_x:\Spec k\to \Spec\sO$ and  $j:\Spec K\to \Spec\sO$, respectively. Then the diagram 
\[
\xymatrix{
\DTM(\sO)\ar[r]^{j^*}\ar[rd]_{i_x^*}&\DTM(K)\ar[d]^{sp_t}\\
&\DTM(k(x))
}
\]
commutes up to natural isomorphism.
\end{prop}

\begin{proof}This follows from Proposition~\ref{prop:Com1} and the functoriality (up to natural isomorphism)  of the equivalence $\DCM(\sN(S))\sim \DTM(S)$.
\end{proof}

\subsection{Compatibility with specialization on motivic cohomology} As above, let $\sO$ be  the local ring of a closed point $x$ on a smooth curve $C$ over $k$, with quotient field $K$ and uniformizing parameter $t$. We have the localization sequence for motivic cohomology
\[
\ldots\to H^n(\sO,\Z(q))\xrightarrow{j^*} H^n(K,\Z(q))\xrightarrow{\partial}H^{n-1}(k(x),\Z(q-1))\xrightarrow{i_{x*}} \ldots
\]
In addition, the parameter $t$ determines the element $[t]\in H^1(K,\Z(1))$.
One defines the {\em specialization homomorphism}
\[
\sp_t:H^n(K,\Z(q))\to H^n(k(x),\Z(q))
\]
by the formula
\[
\sp_t(\alpha):=\partial([t]\cup\alpha).
\]

On the other hand, if $k(x)=k$,  we have the specialization functor
\[
sp_t:\DTM(K)\to \DTM(k(x))
\]
and the natural identifications
\begin{align*}
&H^n(K,\Q(q))\cong \Hom_{\DTM(K)}(\Q,\Q(q)[n])\\
&H^n(k,\Q(q))\cong \Hom_{\DTM(k)}(\Q,\Q(q)[n]).
\end{align*}
Thus the functor $sp_t$ induces the homomorphism
\[
sp_t:\Hom_{\DTM(K)}(\Q,\Q(q)[n])\to \Hom_{\DTM(k)}(\Q,\Q(q)[n])
\]
and hence a new homomorphism
\[
sp_t':H^n(K,\Q(q))\to H^n(k,\Q(q)).
\]

\begin{prop}\label{prop:Compat} $sp_t'$ agrees with the $\Q$-extension of $\sp_t$.
\end{prop}

\begin{proof} Using the equivalence $\DTM(K)\sim\DCM(\sN(K))$ and the canonical identifications
\[
\Hom_{\DCM(K)}(\Q,\Q(q)[n])\cong H^n(\sN(K))\cong \oplus_{q\ge0}H^n(K,\Q(q))
\]
(and similarly for $k$) we need to show that the $\Q$-linear extension of
 $\sp_t$ agrees with the map
 \[
 H^n(sp_t): H^n(\sN(K))\to  H^n(\sN(k))
 \]
 induced by $sp_t:\sN(K)\to \sN(k)$.
 
 For this, take an element $\alpha\in H^n(K,\Z(q))$ and set
 \[
 \bar{\beta}:=\partial\alpha\in H^{n-1}(k,\Z(q-1)).
 \]
 Since $i_x:x\to \Spec\sO$ is split by the structure morphism $\pi:\Spec\sO\to\Spec k$, we can lift $\bar{\beta}$ to $\beta:\pi^*(\bar{\beta})\in H^{n-1}(\sO,\Z(q-1))$. Then
 \[
 \partial([t]\cup\beta)=\partial([t])\cup i_x^*\beta=\bar{\beta},
 \]
 the first identity following from the Leibniz rule and the second from the fact that $\partial([t])=1\in H^0(k,\Z(0))$. Thus
 \[
 \partial(\alpha-[t]\cup\beta)=0,
 \]
 hence there is a class $\gamma\in H^n(\sO,\Z(q))$ with
 \[
 j^*\gamma=\alpha-[t]\cup\beta.
 \]
 We consider $\gamma$ as an element of $H^n(\sN(\sO))$.
 
 By Proposition~\ref{prop:Com1}, we have
 \[
 H^n(i_x^*)(\gamma)=H^n(sp_t)(\alpha-[t]\cup\beta).
 \]

  By the functoriality of the identification 
  \[
  H^n(\sN(-))\cong \oplus_{q\ge0} \Hom_{\DCM}(\sN(-))(\Q,\Q(q))
  \]
   and Proposition~\ref{prop:Com1} it follows that
 \[
\sp_t(j^*\gamma)= H^n(i_x^*)(\gamma)=H^n(sp_t)(j^*\gamma)
\]
so we reduce to showing
\[
\sp_t([t]\cup\beta)=0= H^n(sp_t)([t]\cup\beta).
\]
The first identity follows from $[t]\cup[t]=0$ in $H^2(K,\Q(2))$. For the second, because $sp_t$ is  a morphism in $\sH(\cdga_q)$,  the map $H^*(sp_t)$ is multiplicative, hence it suffices to show that $H^1(sp_t)([t])=0$. 

For this, it follows from the constuction of the map $sp_t:\sN(K)\to \sN(k(x))$ in $\sH(\cdga_\Q)$ that $sp_t$ is natural with respect to Nisnevic neighborhoods $f:(C',x')\to (C,x)$ of $x$, i.e.,
\[
sp_{f^*(t)}\circ f^*=sp_t.
\]
Now, the map $t:(C,x)\to (\A^1,0)$ is clearly a Nisnevic neighborhood of 0 (after shrinking $C$ if necessary) and 
\[
[t]=t^*([T])
\]
where $\A^1=\Spec k[T]$. Thus, we may assume that $C=\A^1$ and $t=T$. But then $[T]$ is a well-defined element of $H^1(\sN(\G_m))$ hence
\[
H^1(sp_t)([T])=i_1^*([T])=[1]=0
\]
by definition of $sp_t:\sN(\sO_{\A^1,0})\to\sN(k)$. This completes the proof.
\end{proof}

\begin{rem} Since $sp'_t$ is multiplicative, as we have already remarked, Proposition~\ref{prop:Compat} gives a sneaky re-proof of the multiplicativity of the specialization homomorphism $\sp_t$
\end{rem}

\subsection{Tangential base-points} As shown in \cite{LevineMTM}, the category $\DTM(k)$ carries a canonical exact {\em weight filtration}. For an Adams-graded \cdga\ $\sA$, the derived category of cell modules $\DCM(\sA)$ carries a natural weight filtration as well; the equivalence $\DCM(\sN(k))\sim \DTM(k)$ given by Spitzweck's representation theorem is compatible with the weight filtrations \cite[Theorem 5.24]{LevineHB}. 

If $\sA$ is cohomologically connected ($H^n(\sA)=0$ for $n<0$ and $H^0(\sA)=\Q\cdot\id$), then $\DCM(\sA)$ carries a $t$-structure, natural among cohomologically connected $\sA$. Finally, if $\sA$ is {\em 1-minimal} then $\DCM(\sA)$ is equivalent to the derived category of the heart of this $t$-structure (see \cite[\S5]{LevineHB}).

Thus, if $\sN(F)$ is   cohomologically connected, then  $\DTM(F)$ has a $t$-structure; the  heart is called the category of mixed Tate motives over $F$, denoted $\MTM(F)$. In fact, $\MTM(F)$ is a  Tannakian category , with natural fiber functor given by the weight filtration; let $ \Gal_\mu(F)$ denote the pro-algebraic group scheme over $\Q$ associated to $\MTM(F)$ by the Tannakian formalism. If $\sN(F)$ is 1-minimal, then  $\DTM(F)$ is equivalent to $D^b(\MTM(F))$, but we won't be using this.

Now let $x$ be a $k$-point on a smooth curve $C$ over $k$, and $t$ a parameter in $\sO_{C,x}$. The specialization functor
\[
sp_t:\DTM(k(C))\to \DTM(k(x))
\]
arises from the map $sp_t:\sN(k(C))\to\sN(k(x))$ in $\sH(\cdga_\Q)$, hence $sp_t$ is compatible with the weight filtrations. When $\sN(k(C))$ and $\sN(k(x))$ are cohomologically connected, 
 $sp_t$ is compatible with the $t$-structures, hence induces an exact functor of Tannakian categories
 \[
 sp_t:\MTM(k(C))\to\MTM(k(x)).
\]
By Tannakian duality, $sp_t$ is equivalent to a homomorphism
\[
\frac{\del}{\del t}_*:\Gal_\mu(k(x))\to \Gal_\mu(k(C)),
\]
called the {\em tangential base-point} associated to the parameter $t$.

\end{document}